\documentclass[11pt]{article}
\usepackage{amsmath, amsthm, amssymb}
\begin{document}
\title{On Asymptotic {\WP} Geometry of Teichm\"{u}ller Space 
of Riemann Surfaces}

\author{Zheng Huang}
\newtheorem{theorem}{Theorem}[section]
\newtheorem*{thm1-10}{Theorem 1 of {\cite {H}}}
\newtheorem*{thm1}{Theorem 1.1}
\newtheorem*{thm2}{Theorem 1.2}
\newtheorem*{thm3}{Theorem 1.3}
\newtheorem*{thm4}{Theorem 1.4}
\newtheorem*{thm5}{Theorem 1.5}

\newtheorem{cor}[theorem]{Corollary}
\newtheorem{lem}[theorem]{Lemma}
\newtheorem{pro}[theorem]{Proposition}
\newtheorem{rem}[theorem]{Remark}

\newcommand{\WP}{Weil-Petersson}
\newcommand{\TS}{Teichm\"{u}ller space}
\newcommand{\Tt}{Teichm\"{u}ller theory}
\newcommand{\hq}{holomorphic quadratic}
\newcommand{\RS}{Riemann surface}
\newcommand{\hm}{harmonic map}
\newcommand{\Sc}{sectional curvature}
\newcommand{\cd}{complex dimension}
\newcommand{\Bd}{Beltrami differential}
\newcommand{\im}{identity map}
\newcommand{\Hd}{Hopf differential}
\newcommand{\ts}{tangent space}
\newcommand{\de}{differential equation}
\newcommand{\mc}{model case}
\newcommand{\tp}{tangent plane}

\maketitle
\tableofcontents 
\section {Introduction}

A remarkable property of hyperbolic surfaces is that many distinct 
complex structures, or equivalently hyperbolic structures, can be 
introduced on a surface. The moduli problem of Riemann asks how many 
distinct complex structures can exist on a closed surface, and this 
problem has been developed into the modern theory of {\TS}.

In this paper, we assume $\Sigma$ is a smooth, closed Riemann 
surface of genus $g$, with $n$ punctures and $3g - 3 + n > 1$. {\TS} 
${\mathcal {T}}_{g,n}$ is the space of hyperbolic metrics (with 
constant curvature $-1$) on $\Sigma$, where two hyperbolic metrics 
$\sigma$ and $\rho$ are equivalent if there is a biholomorphic map 
between $(\Sigma,\sigma)$ and $(\Sigma,\rho)$ in the homotopy class 
of the {\im}.

{\TS} ${\mathcal {T}}_{g,n}$ has its own natural complex structure 
(\cite {Ah}): ${\mathcal {T}}_{g,n}$ is a complex manifold of {\cd} 
$3g - 3 + n > 1$, the co{\ts} at $\Sigma$ is identified with 
$QD(\Sigma)$, the space of {\hq} differentials; and the {\ts} at 
$\Sigma$ is identified with the space of so-called harmonic {\Bd}s.

The {\WP} metric on {\TS} is naturally defined by duality from the 
$L^2$ inner product on $QD(\Sigma)$. This metric is considered as 
one of the natural metrics on {\TS}. Every {\WP} isometry of {\TS} 
is induced by an element of the extended mapping class group when 
$3g - 3 + n > 1$ and $(g,n) \not= (1,2)$ (\cite {MaW}). And the {\WP} 
metric has many interesting geometric properties: for example, it is 
a Riemannian metric with negative {\Sc} (\cite {Wp86} \cite {T87}); 
yet it is incomplete, since not every geodesic can be extended 
indefinitely, the surface $\Sigma$ developes a node when a 
geodesic cannot be further extended (\cite {Mas} \cite {Chu} 
\cite {Wp75}); this metric is K\"{a}hler (\cite {Ah}), and there 
is a negative upper bound $\frac {-1}{2 \pi (g-1)}$, which only 
depends on the topology of the surface, for the holomorphic {\Sc} 
and Ricci curvature (\cite {Wp86}); however, there are no negative 
upper bounds for the {\Sc} (\cite{H}).

Our goal in this paper is to investigate the asymptotic geometry of 
the {\WP} metric on {\TS} and give an estimate on the upper and lower 
bounds for the {\WP} {\Sc} at any point in {\TS}, purely in terms of 
the length of the shortest geodesic on the surface: 
\begin{theorem}
Let $l$ be the length of the shortest geodesic on closed surface 
$\Sigma$, and $K$ be the {\WP} sectional curvature of 
{\TS} $\mathcal {T}$, assuming $dim_{C}{\mathcal {T}} > 1$, 
there exists a constant $C > 0$ such that 
\begin{center}
$-(Cl)^{-1} \le K \le -Cl$.
\end{center}
Moreover, there are {\tp}s with the {\WP} curvatures of 
the orders $O(l)$ and comparable to $l^{-1}$, and hence the {\WP} {\Sc} 
has neither negative upper bound, nor lower bound. 
\end{theorem}

We will firstly prove a result on asymptotic holomorphic {\Sc}. 
Our estimate indicates that the holomorphic {\Sc} tends to 
negative infinity when a core geodesic on the surface is performed 
infinitesimal twists and length shrinking. More specifically, we 
show that
\begin{theorem} 
If the {\cd} of {\TS} $\mathcal {T}$ is greater than 1, then there 
is no negative lower bound for the holomorphic sectional curvature 
of the {\WP} metric. Moreover, let $l$ be the length of the shortest 
geodesic along a path to the frontier of {\TS}, then there 
exists a sequence of {\tp}s with {\WP} holomorphic {\Sc} of the 
order comparable to $l^{-1}$. 
\end{theorem}

We realized that one of the difficulties in estimating curvatures 
is working with the operator $D = -2(\Delta - 2)^{-1}$, which 
appears in Tromba-Wolpert's curvature formula. As pointed out in 
(\cite {Wf89}), there is a natural correspondence of the operator 
$D$ and local variations of the energy of a {\hm} between surfaces. 
By investigating {\hm}s from a nearly noded surface to nearby 
hyperbolic structures, in {\cite {H}}, we showed that even though 
the {\Sc}s are negative, they are not staying away from zero. More 
specifically, we detected an asymptotically flat {\tp}, denoted by 
$\Omega_l$, spanned by two {\Bd}s $\dot{\mu}_0$ and $\dot{\mu}_1$, 
resulting from pinching two independent core geodesics on the 
surface. 
\begin{thm1-10} 
If the {\cd} of {\TS} $\mathcal {T}$ is greater than 1, then the 
{\WP} {\Sc} is not pinched from above by any negative constant. 
Moreover, the {\WP} {\Sc} of $\Omega_l$ is of the order $O(l)$. 
\end{thm1-10}

In theorem 1.2, we are detecting a {\tp}, spanned by {\Bd}s 
$\dot{\mu}_0$ and $i\dot{\mu}_0$, whose curvature 
is asymptotically negative infinity. Following a suggestion of Scott 
Wolpert, we consider a family of {\tp}s $\Omega'_{l}$, spanned by 
$i\dot{\mu}_0$ and $\dot{\mu}_1$, and find that
\begin{theorem}
$\Omega'_{l}$ is asymptotically flat, i.e., its {\WP} curvature 
is of the order $O(l)$.
\end{theorem}

In other words, we are detecting another asymptotically flat {\tp}, 
spanned by $i\dot{\mu}_0$ and $\dot{\mu}_1$. This asymptotic 
flatness results from pinching two nonhomotopic closed geodesics 
on the surface while performing infinitesimal twists on one of 
them. Together with theorem 1 of (\cite {H}), this suggests an 
asymptotic product structure of {\WP} metric, as pointed out by 
Wolpert (\cite {Wp02}). Similarly, we also show that:
\begin{theorem}
The plane $\Omega''_{l}$ spanned by {\Bd}s $i\dot{\mu}_0$ and 
$i\dot{\mu}_1$ is asymptotically flat with respect to the {\WP} 
metric, and its {\WP} {\Sc} is of the order $O(l)$.
\end{theorem}

We notice that theorems 1.3, 1.4 and theorem 1 of (\cite {H}) 
are theorems concerning asymptotic flatness, and the path we take 
towards the frontier of {\TS} (see definition of the frontier 
of {\TS} in $\S 2.1$) is to pinch two short independent geodesics 
on the surface. We can also pinch just one simple closed curve on 
the surface to take a path towards the frontier space. There is 
also a phenomenon of asymptotic flatness when a separating 
geodesic on the surface is pinched. When $\Sigma$ is a closed 
surface with genus at least two and $\gamma_0$ is a separating 
short geodesic on $\Sigma$ with length $l$. Let $\gamma_2$ and 
$\gamma_3$ be two closed geodesics on $\Sigma$ which have fixed 
length $l_0 >> l$, and these two geodesics lie on different sides 
of the shrinking curve $\gamma_0$. We can define a {\tp} 
$\Omega'''_{l}$, spanned by {\Bd}s $\dot{\mu}_2$ and 
$\dot{\mu}_3$, where $\dot{\mu}_2$ and $\dot{\mu}_3$ are obtained 
from infinitesimal twists about curves $\gamma_2$ and $\gamma_3$, 
respectively. 
\begin{theorem}
This plane $\Omega'''_{l}$ is asymptotically flat with respect to 
the {\WP} metric, moreover, its {\WP} {\Sc} is of the order $O(l)$.
\end{theorem}

Here is the more detailed content of this paper. We give the 
necessary background in section 2. Section 3 is devoted to 
proving theorem 1.2. The discussion of this purpose is broken 
into subsections: in $\S 3.1$, we study a harmonic mapping 
problem between hyperbolic cylinders. The family of 
rotationally symmetric {\hm}s will be an approximation to the 
actual family of {\hm}s restricted on the pinching 
neighborhood. We will describe so called ``{\mc} one" in 
$\S 3.2$, namely, we pinch one core geodesic of a cylinder 
into a point, and perform infinitesimal twists about this 
geodesic, and study the asymptotic behavior of the {\hm}s 
between cylinders. We will establish the estimates of terms 
in the curvature formula in {\mc} one; in $\S 3.3$, we 
construct families of maps which have small tension, and are 
close to the {\hm}s resulting from pinching and twisting 
process in $\S 3.2$; finally in $\S 3.4$, we prove theorem 
1.2 based on the estimates in $\S 3.2$ and the construction 
in $\S 3.3$. In section 4 we study the asymptotic flatness 
of the {\WP} metric with two nonhomotopic geodesics on the 
surface are shrinking. We describe {\mc} two in $\S 4.1$, 
estimate curvature terms in $\S 4.2$, and prove theorems 1.3 
and 1.4 in $\S 4.3$. The aim of section 5 is to prove 
curvature bounds for the {\WP} metric. In $\S 5.1$, we 
consider the asymptotic flatness when only one separating 
geodesic on the surface is pinched, and prove theorem 1.5; 
We prove theorem 1.1 in $\S 5.2$.

$\bf{Acknowledgements}$ The author expresses his deepest thanks to 
Mike Wolf for his mentorship, continuous encouragement and many 
fruitful discussion. Theorem 1.1 was brought up to the author by 
Scott Wolpert, the author also wants to thank him for many helpful 
discussion. The author also thanks Yair Minsky for suggesting him 
to prove theorem 1.5.

\section{Harmonic Maps and Teichm\"{u}ller Space}
We will give some background in this section.
\subsection{Teichm\"{u}ller Space of Riemann Surfaces}

Recall that $\Sigma$ is a fixed, oriented, smooth surface of genus $g \geq 1$, 
and $n \geq 0$ punctures where $3g - 3 + n > 1$. We denote hyperbolic metrics 
on the surface $\Sigma$ by $\sigma |dz|^2$ and $\rho |dw|^2$, where $z$ and 
$w$ are conformal coordinates on $\Sigma$. On $(\Sigma,\sigma |dz|^2)$, we 
denote
\begin{center}
$\Delta = \frac {4}{\sigma} \frac {\partial^{2}}{\partial z \partial \bar{z}}$, \\
$K(\rho) = - \frac {2}{\rho} \frac {\partial^{2}}{\partial w \partial \bar{w}} 
log\rho $, \\
$K(\sigma) = -\frac {2}{\sigma} \frac {\partial^{2}}{\partial z \partial \bar{z}}
log\sigma $,
\end{center}
where $\Delta$ is the Laplacian, and $K(\rho)$, $K(\sigma)$ are curvatures of the 
metrics $\rho$ and $\sigma$, respectively.

By the uniformization theorem, the set of all similarly oriented hyperbolic 
structures $M_{-1}$ can be identified with the set of all conformal (or complex) 
structures on $\Sigma$ with the given orientation. And {\TS} 
$\mathcal{T}$ is defined to be the quotient space 
\begin{center}
$\mathcal{T}$$ = {M_{-1}} / Diff_{0}(\Sigma)$
\end{center}

The moduli space of {\RS}s admits the Deligne-Mumford compactification 
(\cite {Mu}), and any element of the 
compactification divisor can be thought of as a {\RS} with nodes, a 
connected complex space where points have neighborhoods complex isomorphic to 
either $\{ |z| < \varepsilon \}$ (regular points) or 
$\{ zw = 0; |z|, |w| < \varepsilon \}$ (nodes). 
We can think of noded surfaces arising as elements of the compactification 
divisor through a pinching process: fix a family of simple closed curves on 
the surface $\Sigma$ such that each component of the complement of the curves 
has negative Euler characteristic. Topologically, the noded surface is the 
result of identifying each curve to the node (\cite {Be}).

{\TS} $\mathcal{T}$ is a complex manifold when $3g-3+n > 1$, and the cotangent 
space at a point $\Sigma \in \mathcal{T}$ is the space of holomorphic 
quadratic differentials $\Phi dz^2$ on $\Sigma$ (\cite {Ah}). The 
{\WP} metric on $\mathcal{T}$ is defined on 
$QD(\Sigma) \cong T^{*}_{\Sigma}\mathcal{T}$ by the $L^2$-norm:
\begin{center}
$||\phi||^2 = \int_{\Sigma} \frac {|\phi|^2}{\sigma}dzd\bar{z}$
\end{center}
where $\sigma |dz|^2$ is the hyperbolic metric on $\Sigma$. By duality, 
we obtain a Riemannian metric on the tangent space to $\mathcal{T}$.

{\TS} with {\WP} metric is not complete (\cite {Chu} \cite {Wp75}). The 
incompleteness is caused by pinching of at least one short geodesic on 
the surface. For the surface with genus at least two, we define 
$\bar{\mathcal{T}}$ to be the {\WP} completion of $\mathcal{T}$, and 
denote $\partial \mathcal{T}$ as the frontier set 
$\bar{\mathcal{T}} \backslash {\mathcal{T}}$. Hence as shown in 
(\cite{Mas}), every point in $\partial \mathcal{T}$ represents a noded 
surface, i.e., the frontier set $\partial \mathcal{T}$ consists of a 
union of lower dimensional {\TS}s, each such space consists of 
topologically reduced {\RS}s (noded surfaces), obtained by pinching 
nontrivial geodesics on the surface.

The curvature tensor of the {\WP} metric is given by (\cite {Wp86})
\begin{center}
$R_{\alpha {\bar{\beta}} \gamma {\bar{\delta}}} =  
(\int_{\Sigma}D(\dot{\mu}_{\alpha} \dot{\bar{\mu}}_{\beta}
\dot{\mu}_{\gamma}\dot{\bar{\mu}}_{\delta}dA) + 
(\int_{\Sigma}D(\dot{\mu}_{\alpha} \dot{\bar{\mu}}_{\delta})
\dot{\mu}_{\gamma}\dot{\bar{\mu}}_{\beta}dA)$
\end{center}
where $dA$ is the area element and $\dot{\mu}$'s are infinitesimal 
{\Bd}s. Here the operator $D = -2(\Delta-2)^{-1}$. It is known that the 
operator $D$ is a self-adjoint compact integral operator with a positive 
kernel, and it is the identity on constant functions.

Then the curvature of $\Omega$ is then given by $R/{\Pi}$, where 
(\cite {Wp86})
\begin{center}
$R = R_{0 {\bar{1}} 0 {\bar{1}}} - R_{0 {\bar{1}} 1 {\bar{0}}} - 
R_{1 {\bar{0}} 0 {\bar{1}}} + R_{1 {\bar{0}} 1 {\bar{0}}}$
\end{center}
and 
\begin{eqnarray*}
\Pi & = & 4 <\dot{\mu}_0,\dot{\mu}_0><\dot{\mu}_1,\dot{\mu}_1> - 
2 |<\dot{\mu}_0,\dot{\mu}_1>|^2 - 2 Re(<\dot{\mu}_0,\dot{\mu}_1>)^2 \nonumber\\
& = & 4 <\dot{\mu}_0,\dot{\mu}_0><\dot{\mu}_1,\dot{\mu}_1> - 
4 |<\dot{\mu}_0,\dot{\mu}_1>|^2
\end{eqnarray*}

It is known that the {\WP} {\Sc}, holomorphic {\Sc} and Ricci curvature 
are all negative (\cite{Wp86} \cite {T87}).
\subsection {Harmonic Maps and Local Variations} 

Our main method is to look at families of {\hm}s between degenerating 
{\RS}s. The method of harmonic maps has been intensively studied as 
an important computational tool in understanding the geometry of 
{\TS}. In particular, the second variation of the energy of the 
{\hm} $w = w(\sigma,\rho)$ with respect to the domain structure 
$\sigma$ (or image structure $\rho$) at $\sigma = \rho$ yields the 
{\WP} metric on $\mathcal{T}$ (\cite {T87}, \cite {Wf89}), and one 
can also re-establish Tromba-Wolpert's curvature tensor formula of 
the {\WP} metric from this method (\cite {J}, \cite {Wf89}). 

For a Lipschitz map 
$w: (\Sigma,\sigma |dz|^2) \rightarrow (\Sigma,\rho |dw|^2)$, 
we define the energy density of $w$ at a point to be 
\begin{center}
$e(w;\sigma,\rho)= \frac {\rho (w(z))}{\sigma (z)} |w_{z}|^2 + 
\frac {\rho (w(z))}
{\sigma (z)} |w_{\bar{z}}|^2 $
\end{center}
and the total energy 
$E(w;\sigma,\rho) = \int_{\Sigma} e(w;\sigma,\rho) \sigma dzd\bar{z}$.

A {\hm} is a critical point of the energy fuctional $E(w;\sigma,\rho)$; it 
satisfies the 
Euler-Lagrange equation, namely,
\begin{center} 
$w_{z \bar{z}} + \frac{\rho_w}{\rho} w_z w_{\bar{z}} = 0.$
\end{center}

The Euler-Lagrange equation for the energy is the condition for the vanishing 
of the {\em{tension}}, which is, in local coordinates,
\begin{center} 
$\tau (w) = \Delta w^{\gamma} + 
^{N} \Gamma_{\alpha \beta}^{\gamma}w_i^{\alpha}w_j^{\beta} 
= 0$
\end{center}

It is fundamental(\cite {Al} \cite {ES} \cite {Ha} \cite {SY} \cite {Sa}) that 
given 
$\sigma, \rho$, there exists a unique {\hm} 
$w: (\Sigma,\sigma) \rightarrow (\Sigma,\rho)$ homotopic to the identity of 
$\Sigma$, and this map is in fact a diffeomorphism. Naturally associated to a 
{\hm} $w: (\Sigma,\sigma |dz|^2) \rightarrow (\Sigma,\rho |dw|^2)$ is a 
quadratic differential $\Phi (\sigma, \rho) dz^2$, which is holomorphic with 
respect to the conformal structure of $\sigma$. This association of a 
quadratic differential to a conformal structure then defines a map 
$\Phi: {\mathcal{T}} \rightarrow QD(\Sigma)$ from {\TS} $\mathcal{T}$ to the 
space of {\hq} differentials $QD(\Sigma)$. This map $\Phi$ is a homeomorphism 
(\cite {Wf89}).

Thus we have a {\hq} differential 
$\Phi dz^2 = \rho w_z {\bar{w}}_{z}dz^2$, and evidently
\begin{center}

$\Phi = 0 \Leftrightarrow w$ is conformal $\Leftrightarrow \sigma = \rho.$

\end{center}
where $\sigma = \rho$ means that $(\Sigma,\sigma)$ and $(\Sigma,\rho)$ are 
the same point in {\TS} $\mathcal{T}$. Also note that the map $w$ can be 
extended to surfaces with finitely many punctures.

We define two auxiliary functions as following:
\begin{center}
$\mathcal{H} = \mathcal{H}$$ (\sigma,\rho) = 
\frac {\rho (w(z))}{\sigma (z)} |w_{z}|^2$
\\
$\mathcal{L} = \mathcal{L}$$ (\sigma,\rho) = 
\frac {\rho (w(z))}{\sigma (z)} |w_{\bar{z}}|^2$
\end{center}

The Euler-Lagrange equation gives that 
\begin{center}
$\Delta log {\mathcal{H}} = -2K(\rho){\mathcal{H}} + 2K(\rho){\cal{L}} + 
2K(\sigma)$.
\end{center}

When we restrict ourselves to the situation when $K(\sigma)=K(\rho)=-1$, 
we will have the following facts (\cite {Wf89}):
\begin{itemize}
\item
The energy density is $e = {\mathcal{H}} + {\mathcal{L}}$;
\item
The Jacobian is ${\mathcal{H}} - {\mathcal{L}}$;
\item
${\mathcal{H}} > 0$;
\item
The {\Bd} $\mu = \frac {w_{\bar{z}}}{w_z} = \frac {\bar{\Phi}}{\sigma \mathcal{H}}$;
\item
$\Delta log {\mathcal{H}} = 2{\mathcal{H}} - 2{\mathcal{L}} - 2$, where 
${\mathcal{H}} \not= 0$; and 
$\Delta log {\mathcal{L}} = 2{\mathcal{L}} - 2{\mathcal{H}} - 2$, where 
${\mathcal{H}} \not= 0$.
\end{itemize}
Now we consider a family of {\hm}s $w(t)$ for $t$ small, where $w(0) =$ id, the 
{\im}. Denote by $\Phi(t)$ the family of {\Hd}s determined by $w(t)$. We rewrite 
$\Delta log {\mathcal{H}} = 2{\mathcal{H}} - 2{\mathcal{L}} - 2$ as 
\begin{center}
$\Delta log {\mathcal{H}}(t) = 2{\mathcal{H}}(t) - 
\frac {2|\Phi(t)|^2}{\sigma^2 {\mathcal{H}}(t)} - 2$
\end{center}
The maximum principle will force all the odd order $t$-derivatives of the 
holomorphic energy ${\mathcal{H}}(t)$ to vanish, since the above equation only 
depends on the modulus of $\Phi(t)$ and not on its argument(\cite {Wf89}), and 
${\mathcal{H}}(t)$ is real-analytic in $t$ (\cite {Wf91}).

For the family of {\hm}s $w(t)$, the domain hyperbolic structure is fixed and 
the target metric is changing. Wolf computed the $t$-derivative of various 
geometric quantities associated with this family $w(t)$, and we collect these 
local variational formulas into:
\begin{lem} (\cite {Wf89})
For the above notations, we have
\begin{itemize}
\item
${\mathcal{H}}(t) \geq 1$, and ${\mathcal{H}}(t) \equiv 1 \Leftrightarrow t = 0$;
\item
${\dot{\mathcal{H}}}(t) = \partial / \partial t^{\alpha}|_{0}$
${\mathcal{H}}(t) = 0$;
\item
${\dot{\mu}} = \partial / \partial t^{\alpha}|_{0} \mu(t) = 
\bar{\Phi}_{\alpha}/{\sigma}$;
\item
${\ddot{\mathcal{H}}}(t) = 
\frac {\partial ^2}{\partial t^{\alpha} \partial t^{\bar\beta}}|_{0}
{\mathcal{H}}(t) = D \frac {\Phi_{\alpha} \Phi_{\bar{\beta}}}{\sigma^2}$. 
\end{itemize}
\end{lem}
With this lemma, and recall that $D = -2(\Delta-2)^{-1}$, we obtain a 
partial {\de} about $\ddot{\mathcal{H}}(t)$ which will play an important 
role in the proof of our theorems:
\begin{center}
$(\Delta -2)(\ddot{\mathcal{H}}(t)) = -2 \frac {\Phi_{\alpha} 
\Phi_{\bar{\beta}}}{\sigma^2}$
\end{center}
\section{Holomorphic Sectional Curvature}

In this section, we will focus on estimating the  {\WP} holomorphic 
{\Sc}, and prove theorem 1.2, namely,
\begin{thm2} 
If the {\cd} of {\TS} $\mathcal {T}$ is greater than 1, then there is no 
negative lower bound for the holomorphic sectional curvature of the {\WP} 
metric. Moreover, let $l$ be the length of the shortest geodesic along a 
path to a boundary point in {\TS}, then there exists a sequence of {\tp}s 
with {\WP} holomorphic {\Sc} of the order comparable to $l^{-1}$. 
\end{thm2}

We organize this section into following subsections. In $\S 3.1$, we will 
study {\hm}s between hyperbolic cylinders, i.e., we consider a family of 
pertubations of the {\im} between two cylinders; in $\S 3.2$, we estimate 
terms in the curvature formula in the {\mc} one where the surface is a long 
cylinder; in $\S 3.3$, we construct a family of $C^{2,\alpha}$ maps between 
surfaces and show that the constructed maps have small tension and are close 
to the {\hm}s we obtain from the pinching and twisting process; in $\S 3.4$, 
we adapt the estimates in $\S 3.2$ and $\S 3.3$ to prove theorem 1.2.

For the sake of simplicity of exposition, we assume that our surface has no 
punctures. We will discuss the case when finitely many punctures are allowed 
in $\S 3.5$. 

\subsection{Harmonic Maps between Cylinders}

In this subsection, we consider the asymptotics of {\hm}s, as pertubations of 
the {\im}, between two hyperbolic cylinders. In particular, using the same 
notation in {\cite{H}}, we study the boundary value problem of harmonically 
mapping the cylinder
\begin{center} 
$M = [l^{-1}sin^{-1}(l), \pi l^{-1}-l^{-1}sin^{-1}(l)] \times [0,1] $ 
\end{center}
with boundary identification  
\begin{center}
$[{\frac {sin^{-1}(l)}{l}}, 
{\frac {\pi}{l}}-{\frac {sin^{-1}(l)}{l}}] \times \{ 0 \} = 
[{\frac {sin^{-1}(l)}{l}}, {\frac {\pi}{l}}-{\frac {sin^{-1}(l)}{l}}] \times \{ 1 \}$ 
\end{center}
where the hyperbolic length element on $M$ is $lcsc(lx)|dz|$, to the cylinder 
\begin{center} 
$N = [L^{-1}sin^{-1}(L), \pi L^{-1}-L^{-1}sin^{-1}(L)] \times [0,1] $ 
\end{center}
with boundary identification  
\begin{center}
$[{\frac {sin^{-1}(L)}{L}}, 
{\frac {\pi}{L}}-{\frac {sin^{-1}(L)}{L}}] \times \{ 0 \} = 
[{\frac {sin^{-1}(L)}{L}}, {\frac {\pi}{L}}-{\frac {sin^{-1}(L)}{L}}] \times \{ 1 \}$ 
\end{center}
where the hyperbolic length element on $N$ is $Lcsc(Lu)|dw|$. Here $l$ and $L$ are the 
lengths of the simple closed core geodesics in the corresponding cylinders.

Let $w = u + iv$ be this {\hm} between cylinders $M$ and $N$, where 
\begin{center}
$u(l,L;x,y) = u(l,L;x), v(x,y) = y .$ 
\end{center}
The Euler-Lagrange equation becomes 
\begin{center}
$u'' = Lcot(Lu)(u'^2-1)$.
\end{center}
with boundary conditions $u({\frac {sin^{-1}(l)}{l}}) = {\frac {sin^{-1}(L)}{L}}$ and 
$u({\frac {\pi}{2l}}) = {\frac {\pi}{2L}}$. Note that both $M$ and $N$ admit an 
anti-isometric reflection about the curves $\{ \frac {\pi}{2l} \} \times [0,1]$ and 
$\{\frac {\pi}{2L}\} \times [0,1]$.

Since the quadratic differential 
$\Phi = \rho w_z {\bar{w}}_{z} = {\frac {1}{4}} L^2 csc^{2}(Lu)(u'^2 -1)$ is 
holomorphic in $M$, we have 
\begin{center}
$0 = {\frac {\partial}{\partial \bar{z}}}(\rho^2 w_z {\bar{w}}_{z}) = 
{\frac {\partial}{\partial {x}}}({\frac {1}{8}}L^2 csc^{2}(Lu)(u'^2 -1))$
\end{center}

Therefore $L^2 csc^{2}(Lu)(u'^2 -1) = c_{0}(l,L)$, where $c_{0}(l,L)$ is independent of 
$x$, and $c_{0}(l,l) = 0$ since $u(x)$ is the {\im} when $L=l$.

So we have 
\begin{center}
$u' = \sqrt {1 + c_{0}(l,L)L^{-2} sin^2 (Lu)}$
\end{center}
with boundary conditions $u({\frac {sin^{-1}(l)}{l}}) = {\frac {sin^{-1}(L)}{L}}$ and 
$u({\frac {\pi}{2l}}) = {\frac {\pi}{2L}}$. So the solution to the Euler-Lagrange equation 
can be derived from the equation
\begin{center}
$\int_{L^{-1}sin^{-1}L}^u {\frac {dv}{\sqrt {1 + c_{0}(l,L)L^2 sin^2 (Lv)}}} = 
x - l^{-1}sin^{-1}(l)$
\end{center}
with $c_{0}(l,L)$ chosen such that 
$\int_{L^{-1}sin^{-1}L}^{\frac{\pi}{2L}}{\frac {dv}{\sqrt {1 + c_{0}(l,L)L^2 sin^2 (Lv)}}} 
= {\frac{\pi}{2l}} - l^{-1}sin^{-1}l$.

It is not hard to show that when $l \rightarrow 0$, the solution 
$u(l,L;x)$ converges to a solution $u(L;x)$ to the ``noded" problem, i.e., 
$M = [1,+\infty)$, where we require $u(L;1) = L^{-1}sin^{-1}(L)$ and 
$\lim_{x \rightarrow +\infty} {u(L;x)} = \frac {\pi}{2L}$.

This ``noded" problem has the explicit solution as following (\cite {Wf91})
\begin{center} 
$u(L;x) = L^{-1}sin^{-1}\{ \frac {1-{\frac{(1-L)}{(1+L)}e^{2L(1-x)}}}{1+
{\frac{(1-L)}{(1+L)}e^{2L(1-x)}}} \} $
\end{center}
with the holomorphic energy
\begin{center} 
${\mathcal{H}}_{0}(L;x) = {\frac{L^{2}x^2}{4}}[\frac {1+{\sqrt{\frac{(1-L)}{(1+L)}}}e^{L(1-x)}}
{1-{\sqrt{\frac{(1-L)}{(1+L)}}}e^{L(1-x)}}]^2$
\end{center}

\subsection{Model Case One}
Consider the surface $\Sigma$ which is developing a node, i.e., we are pinching one 
short closed geodesic $\gamma_0$ on $\Sigma$ to a point $p_0$. We denote $M_0$ its 
pinching neighborhood, i.e., a cylinder described as $M$ in $\S 3.1$ centered at 
$\gamma_0$. In our model case one, the surface is this cylinder $M_0$ with core 
geodesic $\gamma_0$.

We define $M(l, \theta)$ be the surface with two of the Fenchel-Nielsen coordinates, 
namely, the hyperbolic length of $\gamma_0$ being $l$ and the twisting angle of 
$\gamma_0$ being $\theta$. We obtain a point $M(l) = M(l,0)$ in {\TS} 
${\mathcal{T}}_{g}$. Note that as $l$ tends to zero, the surface is developing a node. 
Fix $M(l)$, we vary the length of $\gamma_0$ into length $L = L(t)$, where $L(0) = l$. 
Thus we obtain a family of {\hm}s $W_{0}(t): M(l) \rightarrow M(L(t),0)=M(L(t))$. 
The $t$-derivative of $\mu_0(t)$, the {\Bd} of the the family $W_{0}(t)$, at $t = 0$ 
represents a tangent vector, $\dot{\mu}_0$, of {\TS} ${\mathcal{T}}_{g}$ at $M(l)$; 
And a tangent vector $i{\dot{\mu}}_0$ is obtained by performing infinitesimal 
twists $\theta(t)$ on the shrinking curve $\gamma_0$, where $\theta(0) = 0$, since 
${\frac {d}{d\theta}}(e^{i\theta}{\dot{\mu}}_0)|_{\theta=0}=i{\dot{\mu}}_0$. Hence 
these two tangent vectors $\dot{\mu}_0$ and $i\dot{\mu}_0$ will span a two dimensional 
subspace of the {\ts} $T_{M(l)}{\mathcal{T}}_g$ to ${\mathcal{T}}_g$, therefore we 
obtain a family, $\Omega^{*}_l$, of {\tp}s.
\begin{pro}
The {\WP} holomorphic {\Sc}s of $\Omega^{*}_l$ is of the order comparable to $l^{-1}$. 
Hence the {\WP} curvatures tend to negative infinity as $l$ tends to zero.
\end{pro}
It is easy to see that Proposition 3.1 implies our theorem 1.2.

We denote $\phi_0(t)$ as the {\Hd} corresponding to the cylinder map $w_{0}(t)$ 
in $M_0$. Here $w_{0}(t): M_0(l) \rightarrow M_0(L(t))$ is a family of {\hm}s 
described as $w = u(x)+iy$ in $\S 3.1$. In other regions of the surface, we 
abuse our notations a little bit and still use $\phi_0$ to denote the {\Hd} 
corresponding to {\hm} $W_{0}(t)$. We also denote $\mu_0$ as the 
corresponding {\Bd} to $\phi_0$.

We denote $a = a(l) = l^{-1}sin^{-1}(l)$, and 
$b = b(l) = \pi l^{-1}-l^{-1}sin^{-1}(l)$. And in this paper, $A \sim B$ means 
$A/C \le B \le CA$ for some constant $C > 0$.

So in $M_0$, recalling that $c_{0}(l,L(t)) = L^{2}csc^{2}(Lu)(u'^{2} - 1)$ is 
independent of $x$, we can choose $L(t)$ so that 
${\frac{d}{dt}}|_{t=0}c_{0}(t) = {\frac{d}{dL}}|_{L=l} c_{0}(l,L) = 4$. Thus 
we can assume that 
$\dot{\phi_0} = {\frac{d}{dt}}|_{t=0}({\frac{1}{4}}c_{0}(t)) = 1$ in $M_0$. 
We notice here ${\dot{c}}_0={\frac{d}{dt}}|_{t=0}$ is never zero for all 
positive $l$, otherwise, we would have 
$\dot{\bar {w}}_z = 0$ as $\dot{\phi} = \sigma \dot{\bar {w}}_z$ and hence 
$w$ is a constant map by rotational invariance of the map.

Note that most of the mass of $|\phi_0|$ resides in the thin part associated 
to $\gamma_0$, near where the core geodesic $\gamma_0$ is pinched, and that is 
why the estimate in the long cylinder is critical in our calculation.

We will now estimate all terms in the curvature formula in this model case one 
where we consider the surface as a long cylinder $M_0$. In $M_0$, the 
corresponding {\Bd} is
\begin{center}
$\dot{\mu}_0 = {\frac {d}{dt}}|_{t=0}({\frac {w_{\bar{z}}}{w_{z}}}) = 
\dot{\bar{\phi}}_0/{\sigma}$
\end{center}
and 
$|\dot{\mu}_0|^{2}|_{M_0} = |i{\dot{\mu}}_0 |_{M_0}^{2}= 
|\dot{\bar{\phi}}_0/{\sigma}|^{2}|_{M_0} = l^{-4}sin^{4}(lx)$.

The {\WP} holomorphic {\Sc} is given by the quotient 
$-R_{0\bar{0}0\bar{0}}/|\dot{\mu}_0|^4$. 
And we will estimate $1/|\dot{\mu}_0|^4$ and $|R_{0\bar{0}0\bar{0}}|$ in the 
cylinder $M_0$ in the next two lemmas.
\begin{lem}
$1/|\dot{\mu}_0|^4 \sim l^6$.
\end{lem} 
\begin{proof}
To show this lemma, we use $|\dot{\mu}_0|^{2}|_{M_0} = l^{-4}sin^{4}(lx)$, 
and noticing that $a = a(l) = l^{-1}sin^{-1}(l)$, and 
$b = b(l) = \pi l^{-1}-l^{-1}sin^{-1}(l)$, then we have 
\begin{eqnarray*}
<\dot{\mu}_0,\dot{\mu}_0>|_{M_{0}} & = & \int_{M_{0}}|\dot{\mu}_0|^2\sigma dxdy \\
 & = & \int_{0}^{1}\int_{a}^{b}|\dot{\mu}_0|^2 \sigma dxdy \\
 & = & \int_{0}^{1}\int_{a}^{b}l^{-2}sin^{2}lx dxdy \\
 & = & \frac {\pi}{2} l^{-3} - l^{-3}sin^{-1}(l) \\
 & \sim &  l^{-3}
\end{eqnarray*}

Note that $|\dot{\mu}_0|^4 \sim (<\dot{\mu}_0,\dot{\mu}_0>|_{M_{0}})^2 \sim l^{-6}$, 
which completes the proof of Lemma 3.2.
\end{proof}

Now we are left to estimate $|R_{0\bar{0}0\bar{0}}| = 
\int_{\Sigma}D(|\dot{\mu}_0|^2){|\dot{\mu}_0|}^{2}\sigma dxdy$. The desired 
estimate is to establish
\begin{lem}
$\int_{M_0}D(|\dot{\mu}_0|^2){|\dot{\mu}_0|}^{2}dA \sim l^{-7}$
\end{lem}
\begin{proof}
Firstly, from Lemma 2.1, we have
\begin{center}
$D(|\dot{\mu}_0|^2) = 
-2(\Delta-2)^{-1} \frac {|\dot{\phi_0}|^2}{\sigma^2} $
\end{center}

We recall in $M_0$, the {\Hd} $\phi_0$ is corresponding to the cylinder map 
$w_0(t): (M_0, \sigma) \rightarrow (M_0, \rho(t))$, and $|\dot{\phi_0}| = 1$. As 
in $\S 2.2$ and $ \S 3.1$, we write the holomorphic energy 
${\mathcal{H}} = \frac {\rho (w(z))}{\sigma (z)} |w_{z}|^2$, 
therefore we can write $D(|\dot{\mu}_0|^2) = \ddot{\mathcal{H}}$. Then in the 
cylinder $M_0$, we have 
\begin{equation}
(\Delta-2) \ddot{\mathcal{H}} = -2 \frac {|\dot{\phi_0}|^2}{\sigma^2} = 
-2 l^{-4}sin^{4}(lx)
\end{equation}

A maximum principle argument implies that $\ddot{\mathcal{H}}$ is positive. This 
$\ddot{\mathcal{H}}$ converges to the holomorphic energy ${\ddot{{\mathcal{H}}}_0}$ 
in the ``noded" problem (of $\S 3.1$) when $x$ is fixed but sufficiently large in 
$[a,\pi/2l]$. This convergence guarantees that $\ddot{\mathcal{H}}$ is bounded on 
the compacta in $[a,b]$ and so we can assume that 
$A_{1}(l) = {\ddot{\mathcal{H}}}(a) = {\ddot{\mathcal{H}}}(l^{-1}sin^{-1}l) = O(1) > 0$. 
Then ${\ddot{\mathcal{H}}}(x)$ solves the following {\de}
\begin{equation}
(l^{-2}sin^{2}(lx)){\ddot{\mathcal{H}}}'' -2{\ddot{\mathcal{H}}} = -2l^{-4}sin^{4}(lx)   
\end{equation}
with the conditions
\begin{center}
${\ddot{\mathcal{H}}}(l^{-1}sin^{-1}l) = A_{1}(l), {\ddot{\mathcal{H}}}'(\pi/{2l}) = 0$
\end{center}

Recall from $\S 3.1$ that all the odd order $t$-derivatives of 
${\mathcal{H}}(t)$ vanish. Also notice that 
\begin{center}
$J(x) = \frac {sin^{2}(lx)}{2l^4}$
\end{center}
is a particular solution to equation (2). Hence we can check, by 
the method of reduction of the solutions, the general solution to 
equation (2) with the assigned conditons has the form 
\begin{center}
${\ddot{\mathcal{H}}}(l;x) = J(x) + A_{2}cot(lx) + A_{3}(1-lxcot(lx))$
\end{center}
where coefficients $A_{2} = A_{2}(l)$ and $A_{3} = A_{3}(l)$ are 
constants independent of $x$ and we can check, by substituting the 
solution into the assigned conditions, that $A_{2}$ and $A_{3}$ are of 
the order comparable to $l^{-1}$.

We compute the following
\begin{eqnarray*}
{\int_{M_0}}{\ddot{\mathcal{H}}}(x){|\dot{\mu}_0|}^{2}\sigma dxdy & = & 
{\int^{1}_{0}\int^{b}_{a}}{\ddot{\mathcal{H}}}(x)(l^{-2}sin^{2}(lx))dxdy \\
& = & {\int^{1}_{0}\int^{b}_{a}}J(x)l^{-2}sin^{2}(lx)dxdy \\
& + & {\int^{1}_{0}\int^{b}_{a}}A_{2}cot(lx)l^{-2}sin^{2}(lx)dxdy \\
& + & {\int^{1}_{0}\int^{b}_{a}}A_{3}(1-lxcot(lx))l^{-2}sin^{2}(lx)dxdy \\
& \sim & l^{-7} + O(l^{-4}) + O(l^{-4}) \\
& \sim & l^{-7}
\end{eqnarray*}

Therefore 
\begin{eqnarray}
{\int_{M_0}}D(|\dot{\mu}_0|^2){|\dot{\mu}_0|}^{2}\sigma dxdy & = & 
{\int^{b}_{a}}{\ddot{\mathcal{H}}}{|\dot{\mu}_0|}^{2}\sigma dx \nonumber \\ 
& \sim & l^{-7} 
\end{eqnarray}
This completes the proof of Lemma 3.3.
\end{proof}
\subsection{Construction of Maps}
In last subsection, we estimated the terms in the curvature formula in the model 
case one where the surface is a long cylinder. Essentially, in the pinching 
neighborhood $M_0$, we used the cylinder map $w_0$ (rotationally symmetric {\hm}) 
instead of the actual {\hm} $W_0$ during the computation. Now we are in the 
general setting, i.e., the surface $\Sigma$ is developing a node. In this 
subsection, we will construct a family of maps $G_0(t)$ to approximate the {\hm} 
$W_0(t)$, and the essential parts of this family are the {\im} of the surface 
restricted in the non-cylindral part and the cylinder map in the pinching 
neighborhood. We will also show that this constructed family $G_0(t)$ is 
reasonably close to the {\hm}s $W_0(t)$; hence we can use the estimates we 
obtained in the previous subsection to the general situation.

We recall some of the notations from previous subsections. We still set 
$M_0$ to be the pinching neighborhood of the node $p_0$. Also $W_0(t)$ is the 
{\hm} corresponding to pinching $\gamma_0$ in $M_0$ into length $L = L(t)$, 
where $L(0) = l$, and twisting $\gamma_0$ with angle $\theta(t)$, where 
$\theta(0)=0$. Let $w_0(t)$ be cylinder maps in {\mc} one, we want to 
show that $W_0(t)$ is close to $w_0(t)$ in $M_0$ and is close to {\im} 
outside of $M_0$.

An important feature of twisting the core geodesic $\gamma_0$ is that it will 
change the values of the corresponding $\ddot{\mathcal{H}}$ on the boundary 
of the long cylinder $M_0$. However, we recall from last subsection that the 
convergence of $\ddot{\mathcal{H}}$ to $\ddot{\mathcal{H}}_0$ in the noded 
problem guarantees $\ddot{\mathcal{H}}$ is bounded on the boundary of $M_0$. 
Therefore we can still use model case one to include the consideration of 
twisting the core geodesic, since we do not concern the boundary values of 
$\ddot{\mathcal{H}}$, as long as they are bounded. And in the compact region, 
where is far away from the shrinking and twisting curve $\gamma_0$, and 
changes of value $\ddot{\mathcal{H}}$ will be small, as we will see in next 
subsection.

We denote subsets 
$\Sigma_0 = \{p \in \Sigma: dist(p, \partial M_0) > 1\}$, and define the 
$1$-tube of $\partial M_0$ as 
$B(\partial M_0, 1) = \{p \in \Sigma: dist(p, \partial M_0) \leq 1\}$, 
the intersection region between the long cylinder $M_0$ and compact 
region of the surface. We can construct a $C^{2,\alpha}$ map 
$G_0: \Sigma \rightarrow \Sigma$ such that 
\begin{center}
$G_0(p) = \left \{
\begin{array}{cc}
w_0(t)(p), & p \in M_0 \cap \Sigma_0\\
p, & p \in (\Sigma_0 \backslash M_0) \\
g_{t}(p), & p \in B(\partial M_0, 1)  
\end{array} \right. $
\end{center}
where $g_{t}(p)$ in $B(\partial M_0, 1)$ is constructed so that it satisfies:
\begin{itemize}
\item
$g_{t}(p) = p$ for $p \in \partial (\Sigma_0 \backslash (M_0 \cup B(\partial M_0,1)))$, 
and $g_{t}(p) = w_{0}(t)(p)$ for $p \in \partial (M_0 \cap \Sigma_0)$;
\item
$g_{t}$ is the {\im} when $t = 0$;
\item
$g_{t}$ is smooth and the tension of $g_t$ is of the order $O(t)$.
\end{itemize}

We note that $G_0$ consists of 3 parts. It is the cylinder map of $M_0$ deep into 
the cylinder region in $M_0 \cap \Sigma_0$, the {\im} in the compact region 
$\Sigma_0 \backslash M_0$, and a smooth map in the intersection region 
$B(\partial M_0, 1)$. Among three parts of the constructed map $G_{0}(t)$, 
two of them, the {\im} and the cylinder map, are harmonic hence have zero 
tension; thus the tension of $G_0(t)$ is concentrated in $B(\partial M_0, 1)$. 
From $\S 3.2$, for the cylinder map $w_0 = u(x) + iy$, we have 
$u' = \sqrt{1 + c_0(t)L^{-2}sin^{2}Lu}$, where 
$c_0(0) = 0, \dot{c}_0(0) = 4$. Hence for 
$x \in [l^{-1}sin^{-1}(l), l^{-1}sin^{-1}(l) + 1]$,
\begin{center}
$w_{0,z}(x,y) = {\frac {1}{2}}(u'(x) + 1) 
= {\frac {1}{2}}((1 + O(1)c_{0}(t))^{\frac {1}{2}}+1) = 1 + O(1)t + O(t^2)$\\
$|w_{0,z}(x,y) - 1| = O(t) \rightarrow 0, (t \rightarrow 0)$\\
$|w_{0,z{\bar{z}}}(x,y)| = |{\frac {1}{4}}u''(x)| = O(|Lcot(Lu)(u'^2 - 1)|) = O(t)$
\end{center}

Thus we can require that $|g_{t,z}| \leq C_{2}t$ and 
$|g_{t,z\bar{z}}| \leq C_{2}t$, and the constant $C_2=C_{2}(t,l)$ 
is bounded in both $t$ and $l$ since the coefficient of $t$ for 
$c_0(t)$ is bounded for small $t$ and small $l$. With the local 
formula of the tension in $\S 2.2$, we have $\tau(G_0(t))$, the 
tension of $G_0(t)$ is of the order $O(t)$. Note that these 
constructed maps $G_0(t)$ are not necessarily harmonic.

Now we are about to compare the constructed family $G_0(t)$ and the family of 
{\hm}s $W_0(t)$. To do this, we consider the following function 
$Q_0 = cosh(dist(W_0,G_0)) - 1$. 

\begin{lem}
$dist(W_0, G_0) \leq C_{3}t$ in $B(\partial M_0, 1)$, where the constant 
$C_3 = C_{3}(t,l)$ is bounded for small $t$ and $l$.
\end{lem}
\begin{proof}
First, we claim that $Q_0$ is a $C^2$ function. Notice that both the {\hm} 
$W_0(t)$ and the constructed map $G_0(t)$ are the {\im} when $t = 0$, and both 
families vary smoothly without changing homotopy type in $t$ for sufficiently small 
$|t|$ (\cite {EL}). For all $l > 0$, and for any $\varepsilon > 0$, there exists 
a $\delta$ such that for $|t| < \delta$, we have 
$|W_0(t) - W_0(0)| < {\frac {\varepsilon}{2}}$ and 
$|G_0(t) - G_0(0)| < {\frac {\varepsilon}{2}}$. Therefore the triangular inequality implies 
that $|W_0(t) - G_0(t)| < \varepsilon$. Since $l$ is positive, the Collar Theorem 
(\cite {Bu}) implies that the surface has positive injectivity radius $r$ bounded below, 
and we choose our $\varepsilon << r$, then $Q_0$ is well defined and 
smooth.

We follow an argument in \cite {HW}. For any unit $v \in T^{1}(B(\partial M_0, 1))$, 
the map $G_0$ satifies the inequality $|\|dG_0(v)\| - 1| = O(t)$, hence 
$|dG_0(v)|^2 > 1 - \varepsilon_0$ where $\varepsilon_0 = O(t)$, then we find that for 
any $x \in \Sigma$,
\begin{eqnarray}
\Delta Q_0 & \geq & min\{|dG_0(v)|^2:dG_0(v) \bot \gamma_x \} Q \nonumber \\
 & - & <\tau (G_0), \bigtriangledown d(\bullet,W_0)|_{G_0(x)}>sinh(dist(W_0,G_0))
\end{eqnarray}
where $\gamma_x$ is the geodesic joining $G_0(x)$ to $W_0(x)$ with initial tangent 
vector $-\bigtriangledown d(\bullet,W_0)|_{G_0(x)}$ and terminal tangent vector 
$\bigtriangledown d(G_0(x),\bullet)|_{W_0(x)}$. If $G_0(t)$ does not coincide with 
$W_0(t)$ on $B(\partial M_0, 1)$, we must have all maxima of $Q_0(t)$ on the 
interior of $B(\partial M_0, 1)$, at any such maximum, we apply the inequality 
$|dG_0(v)|^2 > 1 - \varepsilon_0$ to (12) to find
\begin{eqnarray*}
0 \geq \Delta Q_0 & \geq & (1 - \varepsilon_0)Q_0 - 
\tau(G_0)(sinh(dist(W_0, G_0)) 
\end{eqnarray*}

so that at a maximum of $Q_0$, we have 
\begin{center}
$Q_0 \leq {\frac {\tau(G_0)sinhdist(W_0,G_0)}{(1 - \varepsilon_0)}}$
\end{center}

We notice that $Q_0$ is of the order $dist^{2}(W_0,G_0)$ and $sinhdist(W_0,G_0)$ is of the 
order $dist(W_0,G_0)$, this implies that $dist(W_0,G_0)$ is of the order $O(t)$ in 
$B(\partial M_0, 1)$, which completes the proof of Lemma 3.4.
\end{proof}
\begin{rem}
Lemma 3.4 implies that $Q_0(t)$ is of the order $O(t^2)$ in $B(\partial M_0, 1)$. 
\end{rem}

Note that $B(\partial M_0, 1)$ contains the boundary of the cylinder 
$M_0 \cap \Sigma_0$, which we identify with $[a+1,b-1] \times [0,1]$, 
where, again, $a = a(l) = l^{-1}sin^{-1}(l)$, and 
$b = b(l) = \pi l^{-1}-l^{-1}sin^{-1}(l)$. While in the cylinder 
$M_0 \cap \Sigma_0$, we have the inequality 
\begin{eqnarray*}
\Delta Q_0 & \geq & (1 - \varepsilon_0)Q_0 - \tau(G_0)(sinh(dist(W_0, G_0))\\
& = & (1 - \varepsilon_0)Q_0 - \tau(G_0)(tanh(dist(W_0, G_0))(1+Q_0) \\
& = & (1 - \varepsilon_0 - \tau(G_0))Q_0 - \tau(G_0)(tanh(dist(W_0, G_0)) \\
& \geq & 1/2 Q_0 - C_{4}t^2
\end{eqnarray*}
where the constant $C_4$ is bounded for small $t$ and $l$. Therefore we find 
that $Q_0(z,t)$ decays rapidly in $z=(x,y)$ for $x$ close enough to $\pi/2l$. 
Hence we can assume that $dist(W_0,G_0)$ is at most of order $C't$ in $[a+1,b-1]$, 
here $C'=C'(x,l)$ is no greater that $C_{5}x^{-2}$ for $x \in [a+1,\pi/2l]$, and 
no greater than $C_{5}(\pi/l-x)^{-2}$ for $x \in [\pi/2l,b-1]$, where $C_5$ is 
bounded for small $t$ and $l$. Both maps $W_0$ and $G_0$ are harmonic in 
$M_0 \cap \Sigma_0$, so they are also $C^1$ close (\cite {EL}), i.e. we have 
$|W_{0, \bar{z}} - G_{0, \bar{z}}| \leq C_{5}x^{-2}t$ for small $t$ and $l$, when 
$x \in [a+1,\pi/2l]$. Thus we see that 
$|{\dot{W}}_{0,\bar{z}} - {\dot{G}}_{0,\bar{z}}| = C_{5}x^{-3}$, for 
$x \in [a+1,\pi/2l]$. Also, 
$|{\dot{W}}_{0,\bar{z}} - {\dot{G}}_{0,\bar{z}}| = C_{5}(\pi/l-x)^{-3}$, 
for $x \in [\pi/2l,b-1]$.

As before we denote $\phi_0(t)$ as the family of {\Hd}s corresponding to the 
family of {\hm}s $W_{0}(t)$. We also denote $\mu_0$ as the corresponding {\Bd}s in 
$M_0 \cap \Sigma_0$. Write $\phi^{G}_0 = \rho G_{0,z}(t)\bar{G}_{0,z}(t)$. 
Notice that in $M_0 \cap \Sigma_0$, map $G_0$ is the cylinder map hence harmonic, 
so $\phi^{G}_0$ is the {\Hd} corresponding to $G_0$. When $t = 0$ we have 
$W_{0} = G_{0} =$ identity and $\rho = \sigma$, hence we can differentiate 
$\phi^{G}_0 = \rho G_{0,z}(t)\bar{G}_{0,z}(t)$ in $t$ at $t = 0$, and find that 
$|\dot{\phi}^{G}_0 - \dot{\phi}_{W_0}| = 
\sigma |{\dot{W}}_{0,\bar{z}} - {\dot{G}}_{0,\bar{z}}| 
\leq C_{5}l^{2}x^{-3}csc^{2}(lx) = O(1)$ 
for $x \in [a+1,\pi/2l]$, approaching zero as $x$ tends to $\pi/2l$, and 
$|\dot{\phi}^{G}_0 - \dot{\phi}_{W_0}| 
\leq C_{5}l^{2}(\pi/l-x)^{-3}csc^{2}(lx) = O(1)$ 
for $x \in [\pi/2l,b-1]$. Therefore we have proved
\begin{lem}
$|\dot{\phi}^{G}_0 - \dot{\phi}_{W_0}| = O(1)$ for $x \in [a+1,b-1]$.
\end{lem}
\begin{rem}
Recall that our normalization makes $\dot{\phi}_{W_0} = 1$ in $[a,b]$, and 
$\dot{\phi}^{G}_0$ is actually never zero, Lemma 3.6 implies that 
$\dot{\phi}^{G}_0$ is comparable to $\dot{\phi}_{W_0}$.
\end{rem}
\subsection{Proof of Theorem 1.2}

In last subsection, we constructed a family of maps and found 
that these maps are reasonably close to the actual {\hm}s 
between surfaces. Now we are ready to adapt the estimates in 
the {\mc} one to the general setting, and prove the 
Proposition 3.1, which will imply theorem 1.2. 

\begin{proof}[Proof of Proposition 3.1]
We notice that Lemma 3.2 still holds from triangle inequality 
and Lemma 3.6. It will be sufficient to show that Lemma 3.3 
still holds, which implies desired curvature estimate 
immediately.

Now we are about to estimate 
$\int_{\Sigma}D(|\dot{\mu}_0|^2){|\dot{\mu}_0|}^{2}\sigma dxdy$, 
which breaks into 2 integrals as follows:
\begin{eqnarray}
\int_{\Sigma}D(|\dot{\mu}_0|^2){|\dot{\mu}_0|}^{2}dA & = & 
\int_{M_0 \cap \Sigma_0}
D(|\dot{\mu}_0|^2){|\dot{\mu}_0|}^{2}dA \nonumber \\
& + & \int_{\mathcal{K}}D(|\dot{\mu}_0|^2){|\dot{\mu}_0|}^{2}dA
\end{eqnarray}
where $\mathcal{K}$ is the compact set disjoint from $M_0 \cap \Sigma_0$.

For the last integral, from the previous discussion, 
because of the convergence of the {\hm}s to the {\hm}s of 
``noded" problem, we have both $|{\dot{\mu}_0}|$ and 
$|{\dot{\mu}_1}|$ are bounded. The maximum principle 
implies that 
$D(|{\dot{\mu}_0}|^2) = 
{\ddot{\mathcal {H}}} \leq sup \{|{\dot{\mu}_0}|^2 \}$. 
Note that $\mathcal{K}$ is compact, hence we have the second 
integral is of the order of $O(1)$.

From Lemma 2.1, we have $(\Delta-2) \ddot{\mathcal{H}} = 
-2 \frac {|\dot{\phi}_0|^2}{\sigma^2}$, so we rewirte (5) as 
\begin{eqnarray}
\int_{\Sigma}D(|\dot{\mu}_0|^2){|\dot{\mu}_1|}^{2}\sigma dxdy & = & 
\int_{M_0 \cap \Sigma_0}\ddot{\mathcal{H}}{|\dot{\mu}_1|}^{2}
\sigma dxdy + O(1)
\end{eqnarray}

Recall that 
$(\Delta-2) \ddot{\mathcal{H}}^G = -2 \frac {|\dot{\phi}_0^G|^2}{\sigma^2}$, 
where ${\mathcal{H}}^G$ is the holomorphic energy corresponding to the 
{\mc} one when the {\hm} is the cylinder map, and $\phi_0^G$ is the 
quadratic differential corresponding to the constructed map $G_0$. We also 
denote $\mu_0^G$ to be the {\Bd} corresponding to $\phi_0^G$. Also recall, 
from Lemma 3.6, that $|\dot{\phi}_0 - \dot{\phi}_0^G| = O(1)$ in 
$M_0 \cap \Sigma_0$, here $O(1)$ is bounded in $l$ for small $l$. So we can 
set some $\lambda = O(1)$ (bounded in $l$ for small $l$) such that 
$|\dot{\phi}_0|^2 < \lambda^2 |\dot{\phi}_0^G|^2$ and at the boundary of 
$M_0 \cap \Sigma_0$ satisfies 
$\ddot{\mathcal{H}} < \lambda \ddot{\mathcal{H}}^G$. For example, we can 
take 
$\lambda = 1 + max_{\partial {K}}
({\frac {\ddot{\mathcal{H}}}{\ddot{\mathcal{H}}^G}}, 
{\frac {|\dot{\phi}_0|}{|\dot{\phi}_{0}^G|}})$, this $\lambda = O(1)$ 
because at $\partial {\mathcal{K}} = \partial (M_0 \cap \Sigma_0)$, both 
$\ddot{\mathcal{H}}$, $\ddot{\mathcal{H}}^G$, and 
$\frac {|\dot{\phi}_0|}{|\dot{\phi}_{0}^G|}$ are bounded. Therefore, we 
have
\begin{center}
$(\Delta-2) (\ddot{\mathcal{H}} - \lambda \ddot{\mathcal{H}}^G) = 
2 {\frac {\lambda^2|\dot{\phi}_0^G|^2 - |\dot{\phi}_0|^2}{\sigma^2}} > 0$
\end{center}

So $(\ddot{\mathcal{H}} - \lambda \ddot{\mathcal{H}}^G)$ is a subsolution 
to the {\de} $(\Delta-2)Y = 0$, both with bounded boundary conditions. It 
is not hard to see that in the cylindar $M_0$, the solutions to 
$(\Delta-2)Y = 0$ have the form of 
$Y(l;x) = B_{3}cot(lx) + B_{4}(1-lxcot(lx))$, where constants $B_{3}$ and 
$B_{4}$ satisfy that $B_{3} = O(l)$ and $B_{4} = O(l)$. Hence in 
$M_0 \cap \Sigma_0$, we have 
$\ddot{\mathcal{H}} \leq \lambda \ddot{\mathcal{H}}^G + 
B_{3}cot(lx) + B_{4}(1-lxcot(lx))$. 
Now we find that,
\begin{eqnarray*}
\int_{M_0 \cap \Sigma_0}\ddot{\mathcal{H}}{|\dot{\mu}_0|}^2 dA 
& \leq & \int_{M_0 \cap \Sigma_0}(\lambda \ddot{\mathcal{H}}^G 
+ Y(l;x))({|\dot{\mu}_0|}^2)dA  \nonumber \\
& \leq & \int_{M_0 \cap \Sigma_0}\lambda \ddot{\mathcal{H}}^G 
({2|\dot{\mu}_0^G|}^2+2|\dot{\mu}_0-\dot{\mu}_0^G|^2)dA \nonumber \\
& + & \int_{M_0 \cap \Sigma_0}Y(l;x)
({2|\dot{\mu}_0^G|}^2+2|\dot{\mu}_0-\dot{\mu}_0^G|^2)dA 
\end{eqnarray*}

Recalling the computation in $\S 3.3$, we have the following:
\begin{eqnarray*}
\int_{M_0}(\lambda \ddot{\mathcal{H}}^G)|\dot{\mu}_0^G|^2 
\sigma dxdy & = & O(l^{-7})\\
\int_{M_0}Y(l,x)|\dot{\mu}_0^G|^2 \sigma dxdy & = & O(l^{-2})\\
\int_{M_0 \cap \Sigma_0}(\ddot{\mathcal{H}}^G + 
Y(l;x))(|\dot{\mu}_0-\dot{\mu}_0^G|^2)dA & = & O(l^{-7})\\
\end{eqnarray*}

These add up to 
$\int_{\Sigma}D(|\dot{\mu}_0|^2){|\dot{\mu}_0|}^{2}\sigma dxdy = O(l^{-7})$, 
with Lemma 3.2, we have the holomorphic {\Sc} is of the order 
$O(l^6)O(l^{-7}) = O(l^{-1})$.

Noting that 
$(\Delta-2)(\ddot{\mathcal{H}} - Y(l;x)) = -2|\dot{\mu}_0|^2 < 0$, to 
achieve the inequality of the other direction, we find, in 
$M_0 \cap \Sigma_0$, that 
$\ddot{\mathcal{H}} \ge Y(l;x) = B_{3}cot(lx) + B_{4}(1-lxcot(lx))$, 
for some constants $B_{3} \sim l$ and $B_{4} \sim l$. Combining 
with remark 3.7 and above estimates of the integrals, we find that the 
holomorphic {\Sc} is actually of the order comparable to $l^{-1}$. This 
completes the proof of Proposition 3.1. 
Also this proves theorem 1.2.  
\end{proof}
\subsection{Surface with Punctures}
In last four subsections, we proved theorem 1.2. We want to point out here 
that the assumption on the surface $\Sigma$ being compact is not essential. 
In other words, our theorem 1.2, as well as other theorems in this paper, 
are still true when the surface $\Sigma$ has finitely many punctures. Now 
we prove theorem 1.2 in the case of $\Sigma$ being a punctured surface.

The existence of a harmonic diffeomorphism between punctured surfaces has 
been investigated by Wolf (\cite {Wf91}) and Lohkamp (\cite {L}). In 
particular, Lohkamp (\cite {L}) showed that a homeomorphism between 
punctured surfaces is homotopic to a unique harmonic diffeomorphism with 
finite energy, and the {\hq} differential corresponding to the {\hm} in 
the homotopy class of the identity is a bijection between {\TS} of 
punctured surfaces and the space of {\hq} differentials.

In this case, the set $\Sigma \backslash (M_0 \cap \Sigma_0)$ is no longer 
compact. Let ${\mathcal{K}}_0$ be a compact surface with finitely many 
punctures, and $\{{\mathcal{K}}_m\}$ be a compact exhaustion of 
${\mathcal{K}}_0$. We now estimate 
$\int_{{\mathcal{K}}_0}D(|\dot{\mu}_0|^2){|\dot{\mu}_0|}^{2} dA$. Let 
${\mathcal{H}}(t)$ be the holomorphic energy corresponding to the {\hm} 
$w(t): ({\mathcal{K}}_0, \sigma) \rightarrow ({\mathcal{K}}_0, \rho (t))$, 
then ${\mathcal{H}}(t)$ is bounded from above and below, and has nodal 
limit 1 near the punctures (\cite {Wf91}), hence $|\dot{\mu}_0|^2$ has the 
order $o(1)$ near the punctures. To see this, we consider ${\mathcal{K}}_0$ 
as the union of ${\mathcal{K}}_m$ and disjoint union of finitely many 
punctured disks, each equipped hyperbolic metric 
$\frac {|dz|^2} {z^{2}log^{2}z}$. Then 
$|\dot{\mu}_0| = O (|zlog^{2}z|) \rightarrow 0$ as $z$ tends to the p
uncture, since the quadratic differential has a pole of at most the first 
order. We notice that $\partial {\mathcal{K}}_0$ is the boundary of the 
cylinders, where the {\hm}s converge to a solution to the ``noded" problem 
as $l \rightarrow 0$, hence $D(|\dot{\mu}_0|^2) = {\ddot{\mathcal{H}}}(t)$ 
is bounded on $\partial {\mathcal{K}}_0$. Therefore we apply 
Omori-Yau maximum principle (\cite {O} \cite {Y}) to the {\de} 
$(\Delta -2){\ddot{\mathcal{H}}}= -2|\dot{\mu}_0|^2$ on ${\mathcal{K}}_0$ 
and obtain that 
$sup(D(|\dot{\mu}_0|^2)) \leq 
max(sup(|\dot{\mu}_0|^2), max(D(|\dot{\mu}_0|^2))|_{\partial 
{\mathcal{K}}_0}) = O(1)$. 
Hence we have 
\begin{eqnarray*}
\int_{{\mathcal{K}}_0}D(|\dot{\mu}_0|^2){|\dot{\mu}_0|}^{2} dA & \leq & 
\int_{{\mathcal{K}}_0} sup(|\dot{\mu}_0|^2){|\dot{\mu}_0|}^{2} dA \nonumber \\
& \leq & O(1)O(1)Vol({\mathcal{K}}_0) \nonumber \\
& = & O(1)
\end{eqnarray*}
In other words, our proof carries over to the punctured case, which 
completes the proof of theorem 1.2.

We point out that similar argument can apply to the proof of other theorems 
of this paper when finitely many punctures are allowed. For this reason, we 
will assume in the rest of the paper that our surface has no punctures.

\section{Asymptotic Flatness I: Two Curves Pinching}

In this and next sections, we aim to prove theorem 1.1. Firstly we want to 
discuss the asymptotic flatness of the {\WP} metric on {\TS}, when there are  
at most two curves on the surface are shrinking. We deal with the two curves 
pinching case in this section.

As indicated in $\S 2.1$, the frontier space $\partial {\mathcal {T}}$ of 
{\TS} is a union of lower dimensional {\TS}s, each consists of noded 
surfaces obtained by pinching nontrival geodesics on the surface. Hence to 
obtain two tangent vectors near the infinity, at least one geodesic on the 
surface is being pinched. In $\S 4.1$, we describe {\mc} two where the 
surface is a pair of cylinders; in $\S 4.2$, we prove theorem 1.3 in {\mc} 
two; we give general proofs of theorems 1.3 and 1.4 in $\S 4.3$.

\subsection {Model Case Two}
In this subsection, following {\cite {H}}, we discuss the phenomena of 
asymptotic flatness of the {\WP} resulting from pinching two independent 
curves on the surface.

When we pinch two nonhomotopic curves $\gamma_0$ and $\gamma_1$ on 
surface $\Sigma$ to two points, say $p_0$ and $p_1$, the surface 
$\Sigma$ is developing two nodes. We denote $M_0$ and $M_1$ pinching 
neighborhoods for these two geodesics, i.e., two cylinders described 
as $M$ in $\S 3.1$, centered at $\gamma_0$ and $\gamma_1$, respectively.

We define $M(l_0, l_1)$ be the surface with two of the Fenchel-Nielsen 
coordinates, namely, the hyperbolic lengths of $\gamma_0$ and $\gamma_1$, 
are $l_0$ and $l_1$, respectively. When we set the length of these two 
geodesics equal to $l$ simultaneously, we will have a point 
$M(l) = M(l,l)$ in {\TS} ${\mathcal{T}}_{g}$. As $l$ tends to zero, the 
surface is developing two nodes. At this point $M(l)$, there are two 
tangent vectors $i\dot{\mu}_0$ and $\dot{\mu}_1$. Here $\dot{\mu}_0$ is 
the same as described in $\S 3.2$, i.e., we fix $\gamma_1$ in $M_1$ 
having length $l$, and pinch $\gamma_0$ in $M_0$ into length $L = L(t)$, 
where $L(0) = l$. So the $t$-derivative of $\mu_0(t)$ at $t = 0$ 
represents a tangent vector, $\dot{\mu}_0$, of {\TS} ${\mathcal{T}}_{g}$ 
at $M(l)$, and $i\dot{\mu}_0$ is the tangent vector obtained by 
performing infinitesimal twists on the curve $\gamma_0$, as seen 
in $\S 3.2$. We denote the resulting {\hm} by 
$W_{0}(t): M(l,l) \rightarrow M(L(t),l)$. Similarly when we fix 
$\gamma_0$ in $M_0$ having length $l$, and pinch $\gamma_1$ into length 
$L = L(t)$, the $t$-derivative of $\mu_1(t)$ at $t = 0$ represents another 
tangent vector, $\dot{\mu}_1$, at $M(l)$; we denote the resulting {\hm} 
by $W_{1}(t): M(l,l) \rightarrow M(l;L(t))$. We obtain a family of {\tp}s 
$\Omega'_l$, spanned by $i\dot{\mu}_0$ and $\dot{\mu}_1$.

Now we estimate the curvatures of $\Omega'_l$, the result is:
\begin{thm3}
$\Omega'_{l}$ is asymptotically flat, i.e., its {\WP} 
{\Sc} is of the order $O(l)$.
\end{thm3}
Again, as in $\S 3.2$, we will still use rotational symmetric harmonic 
maps to approximate harmonic maps in the cylinder regions. We describe 
our ``{\mc} two": the surface is a pair of hyperbolic cylinders $M_0$ and 
$M_1$, as we are pinching two independent curves $\gamma_0$ and $\gamma_1$.

We denote $\phi_0(t)$ as the Hopf differential corresponding to the 
cylinder map $w_{0}(t)$ in $M_0$, and $\phi_1(t)$ as the Hopf 
differential corresponding to $w_{1}(t)$ in $M_1$. Here 
$w_{0}(t): M_0(l) \rightarrow M_0(L(t))$ and 
$w_{1}(t): M_1(l) \rightarrow M_1(L(t))$ are harmonic maps described as 
$w = u(x)+iy$ in $\S 3.1$. We still use $\phi_0$ to denote the Hopf 
differential corresponding to harmonic map $W_{0}(t)$ in $M_1$, and 
$\phi_1$ as the Hopf differential corresponding to harmonic map 
$W_{1}(t)$ in $M_0$.

We still denote $a = a(l) = l^{-1}sin^{-1}(l)$, and 
$b = b(l) = \pi l^{-1}-l^{-1}sin^{-1}(l)$. 
Recall from $\S 3.2$, we can assume that $|i\dot{\phi_0}| = 1$ in $M_0$, 
and also $\dot{\phi_1}= 1$ in $M_1$. It is not hard to see that 
$|\dot{\phi_1}|{\arrowvert}_{M_0} = \zeta (x,l)$ for $x \in [a,b]$, where 
$\zeta (x,l)$ satisfies that $\zeta (x,l) \leq C_{1}x^{-4}$ for 
$x \in [a,\pi /2l]$, and $\zeta (x,l) \leq C_{1}(\pi/l-x)^{-4}$ for 
$x \in [\pi /2l,b]$, and $\zeta (x,0)$ decays exponentially in 
$[1, +\infty]$. Here $C_{1} = C_{1}(l)$ is positive and bounded as $l$ 
tends to zero. To see this, notice that $\dot{\phi}_1$ is holomorphic and 
$|\dot{\phi}_1|$ is positive, so $log|\dot{\phi}_1|$ is harmonic in the 
cylinder $M_0$. Hence we can express $log|\dot{\phi}_1|$ in a Fourier 
series $\Sigma a_{n}(x)exp(-iny)$, and we compute 
$0 = \Delta log|\dot{\phi}_1| = \Sigma (a''_n - n^{2}a_n)exp(-iny)$. We 
will see later on that $\dot{\phi}_1$ is close to 0 in $M_0$. Hence we 
conclude the properties $\zeta$ has. Similarly, we assume that 
$|i\dot{\phi_0}||_{M_1} = \zeta (x,l)$ for $x \in [a,b]$. We see that most 
of the mass of $|\dot{\phi}_0|$ resides in the thin part associated to $\gamma_0$, 
and most of the mass of $|\dot{\phi}_1|$ resides in the thin part associated to 
$\gamma_1$ (\cite {Mc}).

Therefore the corresponding {\Bd}s are:
\begin{center}
$|i\dot{\mu}_0|^{2}|_{M_0} = 
|\dot{\bar{i\phi_0}}/{\sigma}|^{2}|_{M_0} = l^{-4}sin^{4}(lx)$ \\
$|\dot{\mu}_1|^{2}|_{M_0} = |\dot{\bar{\phi_1}}/{\sigma}|^{2}|_{M_0} 
= l^{-4}sin^{4}(lx)\zeta^2 (x,l) $ \\
$|i\dot{\mu}_0|^{2}|_{M_1} = l^{-4}sin^{4}(lx)\zeta^2 (x,l)$ \\
$|\dot{\mu}_1|^{2}|_{M_1} = l^{-4}sin^{4}(lx) $
\end{center} 

\subsection{Estimates in Model Case Two}
We recall from $\S 2.1$, the curvature tensor is given by 
\begin{center}
$R_{\alpha {\bar{\beta}} \gamma {\bar{\delta}}} =  
(\int_{\Sigma}D(\dot{\mu}_{\alpha} \dot{\bar{\mu}}_{\beta}
\dot{\mu}_{\gamma}\dot{\bar{\mu}}_{\delta}dA) + 
(\int_{\Sigma}D(\dot{\mu}_{\alpha} \dot{\bar{\mu}}_{\delta})
\dot{\mu}_{\gamma}\dot{\bar{\mu}}_{\beta}dA)$
\end{center}
and the curvature of $\Omega'_l$ is $R/{\Pi}$, where
\begin{center}
$R = R_{0 {\bar{1}} 0 {\bar{1}}} - R_{0 {\bar{1}} 1 {\bar{0}}} - 
R_{1 {\bar{0}} 0 {\bar{1}}} + R_{1 {\bar{0}} 1 {\bar{0}}}$ \\
$\Pi = 4 <i\dot{\mu}_0,i\dot{\mu}_0>^{2}<\dot{\mu}_1,\dot{\mu}_1>^{2}-
4|<i\dot{\mu}_0,\dot{\mu}_1>|^2$
\end{center}

Now we estimate terms in the {\mc} two where the surface is a pair of 
cylinders. Similar to Lemma 2 of {\cite {H}}, we have 
\begin{lem}
$1/{\Pi} = O(l^3)$
\end{lem}
\begin{proof}
As in $\S 3.2$, we have that 
$<\dot{\mu}_0,\dot{\mu}_0>|_{M_0} \sim l^{-3}$. And using 
$|\dot{\mu}_1|^{2}|_{M_0} = l^{-4}sin^{4}(lx)\zeta^2 (x,l)$, and 
$\zeta (x,l) \leq C_{1}x^{-4}$ for $x \in [a,\pi /2l]$, we have
\begin{eqnarray*}
<\dot{\mu}_1,\dot{\mu}_1>|_{M_0} & = & 
\int_{M_0}|\dot{\mu}_1|^2 \sigma dxdy \\ 
& = & 2\int_{0}^{1}\int_{a}^{\pi/2l}
l^{-4}sin^{4}(lx)\zeta^2 (x,l) \sigma dxdy \\
& \leq & 2C^{2}_{1}\int_{0}^{1}\int_{a}^{\pi/2l}
l^{-2}sin^{2}(lx)x^{-8} dxdy \\
& = &  O(1) 
\end{eqnarray*}
Also $|<i\dot{\mu}_0,\dot{\mu}_1>| = 
|\int_{\Sigma} i\dot{\mu}_0 \dot{\bar{\mu}}_1 dA|$, hence, 
\begin{eqnarray*}
|<i\dot{\mu}_0,\dot{\mu}_1>|_{M_0} & = & 
|\int_{M_0}i\dot{\mu}_0\dot{\mu}_1 \sigma dxdy|\\
& \leq & C_{1}\int_{M_0}l^{-2}sin^{2}(lx)x^{-4} dxdy\\
& = & O(1)
\end{eqnarray*}
Note that ${\Pi} \geq (4<i\dot{\mu}_0,i\dot{\mu}_0><\dot{\mu}_1,\dot{\mu}_1> - 
4(|<i\dot{\mu}_0,\dot{\mu}_1>|))^2|_{M_0} \sim l^{-3}$, which completes 
the proof of Lemma 4.1.
\end{proof}

From Lemma 4.1, we have 
\begin{equation*}
|R|/{\Pi} = O(|R|/(l^{-3})) = O(|R|l^3)
\end{equation*}
Now we are left to estimate $|R|$. Similar to Lemma 3 of {\cite {H}}, we have
\begin{lem}
$|R| \leq 8 \int_{\Sigma}D(|i\dot{\mu}_0|^2){|\dot{\mu}_1|}^{2}\sigma dxdy$ 
\end{lem}
\begin{proof}
Note that $D = -2(\Delta - 2)^{-1}$ is self-adjoint, hence we have 
\begin{center}
${\int_{\Sigma}}D(|i\dot{\mu}_0|^2){|\dot{\mu}_1|}^{2}\sigma dxdy = 
{\int_{\Sigma}}D(|\dot{\mu}_1|^2){|i\dot{\mu}_0|}^{2}\sigma dxdy$
\end{center}
Therefore,
\begin{eqnarray*}
R & = & R_{0 {\bar{1}} 0 {\bar{1}}} - R_{0 {\bar{1}} 1 {\bar{0}}} - 
R_{1 {\bar{0}} 0 {\bar{1}}} + R_{1 {\bar{0}} 1 {\bar{0}}}\\
& = & 
2 \int_{\Sigma}D(i\dot{\mu}_0\dot{\bar{\mu}}_1)
i\dot{\mu}_0\dot{\bar{\mu}}_1 \sigma dxdy + 
2\int_{\Sigma}D(\dot{\mu}_{1}(-i\dot{\bar{\mu}}_0))
\dot{\mu}_{1}(-i\dot{\bar{\mu}}_0)\sigma dxdy\\
& - & \int_{\Sigma}D(|i\dot{\mu}_0|^2){|\dot{\mu}_1|}^{2}\sigma dxdy - 
\int_{\Sigma}D(|\dot{\mu}_1|^2){|i\dot{\mu}_0|}^{2}\sigma dxdy\\
& - & \int_{\Sigma}D(i\dot{\mu}_{0}\dot{\bar{\mu}}_1)
\dot{\mu}_{1}(-i\dot{\bar{\mu}}_0)\sigma dxdy - 
 \int_{\Sigma}D(\dot{\mu}_{1}(-i\dot{\bar{\mu}}_0))
i\dot{\mu}_{0}\dot{\bar{\mu}}_{1}\sigma dxdy\\
& = & -6 \int_{\Sigma}D(\dot{\mu}_{0}\dot{\mu}_1)
\dot{\mu}_{1}\dot{\mu}_0 \sigma dxdy - 
2 \int_{\Sigma}D(|i\dot{\mu}_0|^2){|\dot{\mu}_1|}^{2}\sigma dxdy  
\end{eqnarray*}
The last equality follows from that here $\dot{\mu_0}$ and $\dot{\mu_1}$ 
are real functions. And we have $|D(\dot{\mu}_{0}\dot{\mu}_1)| \leq 
|D(|\dot{\mu}_0|^2)|^{1/2}|D(|\dot{\mu}_1|^2)|^{1/2}$ from lemma 4.3 of 
{\cite {Wp86}}. This proves Lemma 4.2. 
\end{proof}

From Lemmas 4.1 and 4.2, we will need to establish next estimate to prove 
theorem 1.3: 
\begin{lem}
$\int_{\Sigma}D(|i\dot{\mu_0}|^2)
|\dot{\mu}_1|^{2}\sigma dxdy = O(l^{-2})$
\end{lem}
\begin{proof}
The surface is a pair of two cylinders, so we need to estimates two 
integrals: $\int_{M_0}D(|i\dot{\mu}_0|^2)|\dot{\mu}_1|^{2}dA$ and 
$\int_{M_1}D(|i\dot{\mu}_0|^2)|\dot{\mu}_1|^{2}dA$.

For the first integral 
$\int_{M_0}D(|i\dot{\mu_0}|^2)|\dot{\mu}_1|^{2}dA$, we write 
$D(|\dot{\mu}_0|^2) = \ddot{\mathcal{H}}$, where ${\mathcal{H}}(t)$ 
is the holomorphic energy for the family of {\hm}s 
$w_{0}(t)$, as in $\S 3.2$. Then in $M_0$, we have 
$(\Delta-2) \ddot{\mathcal{H}} = -2\frac {|\dot{\phi}_0|^2}{\sigma^2}
= -2l^{-4}sin^{4}(lx)$. We recall that $\ddot{\mathcal{H}}$ is bounded 
on the compacta in $[a,b]$ and so we can assume that 
$A_{1}(l) = {\ddot{\mathcal{H}}}(a) = 
{\ddot{\mathcal{H}}}(l^{-1}sin^{-1}l) = O(1) > 0$. 
Then ${\ddot{\mathcal{H}}}(x)$ solves the {\de} 
\begin{equation*}
(l^{-2}sin^{2}(lx))
{\ddot{\mathcal{H}}}'' -2{\ddot{\mathcal{H}}} = -2l^{-4}sin^{4}(lx)   
\end{equation*}
with the conditions
\begin{center}
${\ddot{\mathcal{H}}}(l^{-1}sin^{-1}l) = A_{1}(l), 
{\ddot{\mathcal{H}}}'(\pi/{2l}) = 0$
\end{center}
We remark here that infinitesimal twists only change the boundary 
values for $\ddot{\mathcal{H}}$, which will stay bounded.

We have solved above {\de} in $\S 3.2$, and the general solutions have 
the following form
\begin{center}
${\ddot{\mathcal{H}}}(l;x) = 
\frac {sin^{2}(lx)}{2l^4} + A_{2}cot(lx) + A_{3}(1-lxcot(lx))$
\end{center}
where coefficients $A_{2} = A_{2}(l)$ and $A_{3} = A_{3}(l)$ are 
constants that satisfy 
\begin{center}
$A_{2} = \frac {\pi}{2}A_{3} = O(l^{-1})$
\end{center}
Noticing that $\zeta(x,l) \leq C_{1}x^{-4}$ in $[a,\pi/2l]$, we compute 
the first integral:
\begin{eqnarray}
{\int_{M_0}}D(|\dot{\mu}_0|^2){|\dot{\mu}_1|}^{2}\sigma dxdy & = & 
{\int_{M_0}}{\ddot{\mathcal{H}}}(x){|\dot{\mu}_1|}^{2}\sigma dxdy \nonumber \\
& = & {\int^{1}_{0}\int^{b}_{a}}{\ddot{\mathcal{H}}}(x)
(l^{-2}sin^{2}(lx))\zeta^{2}(x,l)dxdy \nonumber \\
& \leq & 2C^{2}_{1}(\int^{\pi/2l}_{a}
{\frac {sin^{2}(lx)}{2l^4}}l^{-2}sin^{2}(lx)x^{-8}dx \nonumber \\
& + & {\int^{\pi/2l}_{a}}A_{2}cot(lx)l^{-2}sin^{2}(lx)x^{-8}dx \nonumber \\
& + & {\int^{\pi/2l}_{a}}A_{3}(1-lxcot(lx))l^{-2}sin^{2}(lx)x^{-8}dx) \nonumber \\
& = & O(l^{-2}) + O(l^{-2}) + O(l^{-2}) \nonumber \\
& = & O(l^{-2})
\end{eqnarray}
Now let us look at the second integral 
$\int_{M_1}D(|\dot{\mu}_0|^{2})|\dot{\mu}_1|^{2}\sigma dxdy$. 
Note that in $M_1$, which we identify with $[a,b] \times [0,1]$, we have 
$D(|\dot{\mu}_0|^{2}) = \ddot{\mathcal{H}}$, here $\mu_0(t)$ and 
$\mathcal{H}$ come from the harmonic map 
$W_0(t): M(l,l) \rightarrow M(L(t),l)$, where 
$|\dot{\mu}_0||_{M_1} = l^{-2}sin^{2}(lx)\zeta (x,l)$, and 
${\mathcal{H}}(t)$ solves 
\begin{equation}
(l^{-2}sin^{2}(lx))\ddot{\mathcal{H}}'' -2\ddot{\mathcal{H}} = 
-2l^{-4}sin^{4}(lx)\zeta^{2}(x,l)
\end{equation}
with the conditions
\begin{center}
$\ddot{\mathcal{H}}(l^{-1}sin^{-1}l) = B_1(l)$, 
$\ddot{\mathcal{H}}(\pi/{l}-l^{-1}sin^{-1}l) = B_2(l)$.
\end{center}
Here $B_1(l)$ and $B_2(l)$ are positive and bounded as $l$ tends to 
zero, since $\ddot{\mathcal{H}}$ converges to the holomorphic energy 
in the ``noded" problem (of $\S 3.2$, when $M_1 = [1,\infty)$). We 
recall that $|\dot{\phi}_1||_{M_0} = \zeta (x,l)$, where 
$\zeta (x,l) \leq C_{1}x^{-4}$ for $x \in [a,\pi/2l]$ 
and $\zeta (x,l) \leq C_{1}(\pi/l-x)^{-4}$ for $x \in [\pi/2l,b]$.

Consider the equation 
\begin{equation}
(l^{-2}sin^{2}(lx))Y'' -2Y = 0 
\end{equation}
with the boundary conditions that satisfy 
\begin{center}
$Y(l^{-1}sin^{-1}l) = O(1)$ and $Y(\pi/{l}-l^{-1}sin^{-1}l) = O(1)$, 
as $l \rightarrow 0$.
\end{center}
We choose $h(x) = -x^{-2}$ and claim that $\ddot{\mathcal{H}} - h$ is 
a subsolution to (9) for $x \in [a,b]$. To see this, noticing that 
$2|\dot {\phi}_0|^{2}|_{M_1}l^{-4}sin^{4}(lx)$ decays rapidly as 
$x \rightarrow b$ for small $l$, we have
\begin{equation}
(l^{-2}sin^{2}(lx))(\ddot{\mathcal{H}}-h(x))'' -2(\ddot{\mathcal{H}}-h(x)) 
= ({\frac{6sin^{2}(lx)}{l^{2}x^4}}-{\frac{2}{x^2}}) - 
2 |{\dot\mu}_0|^2> 0
\end{equation}

Notice that $\lambda Y(x)$ solves the equation (9) if $Y(x)$ does, for 
any constant $\lambda$. So up to multiplying by a bounded constant, we 
have $Y|_{\partial M_1} > (\ddot{\mathcal{H}} - h(x))|_{\partial M_1}$. 
Hence $\ddot{\mathcal{H}} - h$ is a subsolution to (9) for 
$x \in [a,b]$, while the solutions to (9) have the form  
\begin{center}
$Y(l;x) = B_{3}cot(lx) + B_{4}(1-lxcot(lx))$
\end{center} 
where constants $B_{3} = B_{3}(l)$ and $B_{4} = B_{4}(l)$ satisfy, from 
the boundary conditions for the equation (9), that
\begin{center}
$B_{3} = O(l)$ and $B_{4} = O(l)$
\end{center}

Therefore in $[a,b]$, We have $\ddot{\mathcal{H}} \leq h(x) + Y(x)$. Now,
\begin{eqnarray*}
{\int_{M_1}}Y(x){|\dot{\mu}_1|}^{2}\sigma dxdy & = & 
{\int^{1}_{0}\int^{b}_{a}}Y(x)(l^{-2}sin^{2}(lx))dxdy \\
& = & {\int^{b}_{a}}B_{3}cot(lx)(l^{-2}sin^{2}(lx))dx \\ 
& + & {\int^{b}_{a}}B_{4}(1-lxcot(lx))(l^{-2}sin^{2}(lx))dx \\
& = & O(l^{-2}) + O(l^{-2}) \\
& = & O(l^{-2})
\end{eqnarray*}

Now we compute the second integral:
\begin{eqnarray}
{\int_{M_1}}D(|i\dot{\mu_0}|^{2}){|\dot{\mu}_1|}^{2}dA & = & 
{\int^{1}_{0}\int^{b}_{a}}
\ddot{\mathcal{H}}(x){|\dot{\mu}_1|}^{2}dA \nonumber\\
& \leq & {\int^{1}_{0}\int^{b}_{a}}(Y(x) + h(x))
{|\dot{\mu}_1|}^{2}\sigma dxdy  \nonumber\\
& = & O(l^{-2}) + O(l^{-1}) = O(l^{-2}). 
\end{eqnarray}
This proves Lemma 4.3. 
\end{proof}
Combined with the estimates of Lemmas 4.1, 4.2 
and 4.3, we proved theorem 1.3 when the surface is a pair of cylinders.
\subsection {Proof of Theorems 1.3 and 1.4}

We are now ready to complete the proof of theorem 1.3. Recall from 
$\S 3.3$, we constructed a family of maps $G_{0}(t)$ to approximate 
$W_0(t)$. The map $G_{0}(t)$ consists of three parts: the cylinder map 
$w_0$ of $M_0$ deep into the cylinder region in $M_0 \cap \Sigma_0$, 
the {\im} in the compact region $\Sigma_0 \backslash M_0$, and a 
smooth map interplaying with the {\im} and $w_0$ in the intersection 
region $B(\partial M_0, 1)$. Similarly, we obtain a family of maps 
$G_{1}(t)$ such that it is $w_1$ in $M_1 \cap \Sigma_1$, the {\im} in 
$\Sigma_1 \backslash M_1$, and a smooth map in $B(\partial M_1, 1)$. 
Here $\Sigma_0 = \{p \in \Sigma: dist(p, \partial M_0) > 1\}$, and 
$B(\partial M_0, 1) = \{p \in \Sigma: dist(p, \partial M_0) \leq 1\}$.

\begin{proof}[Proof of Theorem 1.3]
Parallel to the discussion in $\S 3.3$, we know that both families 
$G_{0}(t)$ and $G_{1}(t)$ are close to the families of {\hm}s $W_0(t)$ 
and $W_1(t)$, respectively (Lemma 3.4). Also
$|\dot{\phi}_{G_0} - \dot{\phi}_{W_0}| = O(1)$ and 
$|\dot{\phi}_{G_1} - \dot{\phi}_{W_1}| = O(1)$ both hold for 
$x \in [a+1,b-1]$ (Lemma 3.6).

It is not hard to see that Lemma 4.1 still holds, so we need to establish 
the same estimate as in Lemma 4.3 in order to prove theorem 1.3 in the 
general case, i.e., we need to estimate $\int_{\Sigma}D(|i\dot{\mu_0}|^2)
|\dot{\mu}_1|^{2}\sigma dxdy$, which breaks into three integrals as 
following:
\begin{eqnarray}
\int_{\Sigma}D(|i\dot{\mu_0}|^2){|\dot{\mu}_1|}^{2}\sigma dxdy & = & 
\int_{M_0 \cap \Sigma_0}D(|i\dot{\mu_0}|^2)
{|\dot{\mu}_1|}^{2}\sigma dxdy \nonumber \\
& + & 
\int_{M_1 \cap \Sigma_1}D(|i\dot{\mu_0}|^2)
{|\dot{\mu}_1|}^{2}\sigma dxdy \nonumber\\ 
& + & \int_{{\mathcal{K}}}D(|i\dot{\mu_0}|^2)
{|\dot{\mu}_1|}^{2}\sigma dxdy \nonumber \\ 
& = & \int_{M_0 \cap \Sigma_0}\ddot{\mathcal{H}}
{|\dot{\mu}_1|}^{2}\sigma dxdy  \nonumber \\
& + & \int_{M_1 \cap \Sigma_1}\ddot{\mathcal{H}}
{|\dot{\mu}_1|}^{2}\sigma dxdy + O(1)
\end{eqnarray}
where ${\mathcal{K}}$ is the compact set disjoint from 
$(M_0 \cap \Sigma_0) \cup (M_1 \cap \Sigma_1)$.

Notice that 
$\partial {\mathcal{K}} = \partial 
((M_0 \cap \Sigma_0) \cup (M_1 \cap \Sigma_1))$, 
so parallel to the discussion in $\S 3.4$, the maximum principle will 
force that $\ddot{\mathcal{H}} \leq \lambda \ddot{\mathcal{H}}^G + Y(l,x)$ 
in $(M_0 \cap \Sigma_0) \cup (M_1 \cap \Sigma_1)$, for some bounded 
constant $\lambda$, where $Y(l,x) = B_{3}cot(lx) + B_{4}(1-lxcot(lx))$ 
is the solution to $(\Delta -2)Y=0$, and constants $B_{3}$ and $B_{4}$ 
satisfy that $B_{3} = O(l)$ and $B_{4} = O(l)$. Apply this to (12):
\begin{eqnarray*}
\int_{\Sigma}D(|\dot{\mu}_0|^2){|\dot{\mu}_1|}^2 dA & = & 
\int_{(M_0 \cap \Sigma_0) \cup (M_1 \cap \Sigma_1)}
\ddot{\mathcal{H}}{|\dot{\mu}_1|}^2 dA  + O(1) \nonumber \\
& \leq & \int_{(M_0 \cap \Sigma_0)}(\lambda \ddot{\mathcal{H}}^G 
+ Y(l;x))({|\dot{\mu}_1|}^2)dA \nonumber \\
 & + & \int_{(M_1 \cap \Sigma_1)}(\lambda \ddot{\mathcal{H}}^G 
+ Y(l;x))({|\dot{\mu}_1|}^2)dA + O(1) \nonumber \\
& \leq & \int_{M_0}(\lambda \ddot{\mathcal{H}}^G 
+ Y(l;x))({2|\dot{\mu}_1|}^2+2|\dot{\mu}_1-\dot{\mu}_1^G|^2)dA \nonumber \\
 & + & \int_{M_1}(\lambda \ddot{\mathcal{H}}^G 
+ Y(l;x))({2|\dot{\mu}_1|}^2+2|\dot{\mu}_1-\dot{\mu}_1^G|^2)dA \nonumber \\
 & + & O(1)
\end{eqnarray*}

From the calculation in $\S 3.3$, and $|\dot{\mu}_1-\dot{\mu}_1^G| = 
|\dot{\phi}_1-\dot{\phi}_1^G|/\sigma$, where 
$|\dot{\phi}_1-\dot{\phi}_1^G| \leq C_{5}l^{2}x^{-2}csc^{2}(lx)$ for 
$x \in [a+1,\pi /2l]$, we have the following:
\begin{eqnarray*}
\int_{M_0}(\lambda \ddot{\mathcal{H}}^G)
|\dot{\mu}_1^G|^2 dA & = & O(l^{-2})\\
\int_{M_1}(\lambda \ddot{\mathcal{H}}^G)
|\dot{\mu}_1^G|^2 \sigma dxdy & = & O(l^{-2})\\
\int_{M_0}Y(l,x)|\dot{\mu}_1^G|^2 \sigma dxdy & = & O(1)\\
\int_{M_1}Y(l,x)|\dot{\mu}_1^G|^2 \sigma dxdy & = & O(l^{-2})\\
\int_{M_0 \cap \Sigma_0}(\ddot{\mathcal{H}}^G + Y(l;x))
(|\dot{\mu}_1-\dot{\mu}_1^G|^2)dA & = & o(l^{-2})\\
\int_{M_1 \cap \Sigma_1}(\ddot{\mathcal{H}}^G + Y(l;x))
(|\dot{\mu}_1-\dot{\mu}_1^G|^2)dA & = & o(l^{-2})\\
\end{eqnarray*}
Applying these estimates to (12), we have
\begin{center}
$\int_{\Sigma}D(|i\dot{\mu_0}|^2){|\dot{\mu}_1|}^{2}dA = O(l^{-2})$.
\end{center}
This and Lemma 4.1 complete the proof of theorem 1.3. 
\end{proof}

To prove theorem 1.4, we consider the plane $\Omega''_{l}$, spanned by 
{\Bd}s $i\dot{\mu}_0$ and $i\dot{\mu}_1$. Again we find that 
Lemma 4.1 still holds, i.e., $1/{\Pi} = O(l^3)$, where 
$\Pi = 4 <i\dot{\mu}_0,i\dot{\mu}_0>^{2}<i\dot{\mu}_1,i\dot{\mu}_1>^{2}-
4|<i\dot{\mu}_0,i\dot{\mu}_1>|^2$. A direct calculation shows that $R$ 
is given by
\begin{eqnarray*}
R & = & R_{0 {\bar{1}} 0 {\bar{1}}} - R_{0 {\bar{1}} 1 {\bar{0}}} - 
R_{1 {\bar{0}} 0 {\bar{1}}} + R_{1 {\bar{0}} 1 {\bar{0}}}\\
& = & 
2 \int_{\Sigma}D(i\dot{\mu}_{0}(-i)\dot{\bar{\mu}}_1)
\dot{\mu}_{0}\dot{\bar{\mu}}_1 dA + 
2\int_{\Sigma}D(i\dot{\mu_1}(-i\dot{\bar{\mu}}_0))
\dot{\mu_1}(\dot{\bar{\mu}}_0)dA\\
& - & \int_{\Sigma}D(|i\dot{\mu_0}|^2){|i\dot{\mu_1}|}^{2}dA - 
\int_{\Sigma}D(|i\dot{\mu_1}|^2){|i\dot{\mu_0}|}^{2}dA\\
& - & \int_{\Sigma}D(i\dot{\mu_0}(-i)\dot{\bar{\mu}}_1)
i\dot{\mu_1}(-i\dot{\bar{\mu}}_0)dA - 
 \int_{\Sigma}D(i\dot{\mu_1}(-i\dot{\bar{\mu}}_0))
i\dot{\mu_0}(-i)\dot{\bar{\mu}}_{1}dA\\
& = & 2 \int_{\Sigma}D(\dot{\mu}_{0}\dot{\mu}_1)
\dot{\mu}_{1}\dot{\mu}_{0}dA - 
2 \int_{\Sigma}D(|i\dot{\mu}_0|^2){|i\dot{\mu}_1|}^{2}dA  
\end{eqnarray*}
Hence we have a result similar to that of Lemma 4.2:
\begin{lem}
$|R| \leq 4 \int_{\Sigma}D(|i\dot{\mu}_0|^2){|i\dot{\mu}_1|}^{2}dA$ 
\end{lem}
Now the curvature of $\Omega'''_{l}$ is $R/{\Pi}$. An argument parallel to 
that of proving theorem 1.3 leads to a proof of theorem 1.4. 
\qed
\section {Asymptotic Flatness II: Curvature Bounds}
The goal of this section is to prove theorem 1.5 and, finally, theorem 1.1. 
We continue to discuss the asymptotic flatness of the {\WP} metric when 
there is only one shrinking geodesic on the surface in $\S 5.1$, then give 
a proof of theorem 1.1 in $\S 5.2$.

In this section, we always use $l$ to denote the length of the shortest 
geodesic on the surface. 

\subsection {One Curve Pinching}
In this subsection, we assume that $\Sigma$ is a closed surface of genus 
$g$ at least two, and $\gamma_0$ is a closed separating short geodesic 
on the surface $\Sigma$ with length $l$. As $l$ tends to zero, the 
surface is developing a separating node. Let $\gamma_2$ and $\gamma_3$ be 
two closed geodesics on the different sides of curve $\gamma_0$. Without 
loss of generality, we assume $l(\gamma_2) = l(\gamma_3) >> l$.

As in previous sections, we denote $M(l,\theta_2, \theta_3)$ as the 
surface with three of the Fenchel-Nelsen coordinates, the hyperbolic 
length of $\gamma_0$, the twisting angles of $\gamma_2$ and $\gamma_3$, 
being $l,\theta_2, \theta_3$, respectively.

We denote $M_0$ as a cylinder centered at $\gamma_0$, while $M_2$ and $M_3$ 
are cylinders centered at $\gamma_2$ and $\gamma_3$, respectively. Note 
that as the length of $\gamma_0$ is kept very short, the cylinder $M_0$ 
becomes very long and similar to the cylinder described as $M$ in $\S 3.1$. 
Now we have two families of {\hm}s: the family 
$W_{2}(t): M(l) \rightarrow M(l, \theta_{2}(t),0)$, which fixes $\gamma_3$, 
keeps $\gamma_0$ very short at length $l$, and twists $\gamma_2$ at angle 
$\theta_{2}(t)$; and the family 
$W_{3}(t): M(l) \rightarrow M(l, \theta_{3}(t),0)$, which fixes $\gamma_2$, 
keeps $\gamma_0$ very short at length $l$, and twists $\gamma_3$ at angle 
$\theta_{3}(t)$. We denote $\dot{\mu}_2$ and $\dot{\mu}_3$ as infinitesimal 
{\Bd}s corresponding to $W_{2}(t)$ and $W_{3}(t)$, respectively. Also 
$\dot{\phi}_2$ and $\dot{\phi}_3$ are corresponding infinitesimal {\Hd}s. 
Now we obtain a {\tp} $\Omega'''_l$ at a point $M(l) = M(l,0,0)$ 
in {\TS}, spanned by tangent vectors $\dot{\mu}_2$ and $\dot{\mu}_3$, and 
we will show that
\begin{thm5}
This plane $\Omega'''_{l}$ is asymptotically flat with respect to the 
{\WP} metric, moreover, its {\WP} {\Sc} is of the order $O(l)$.
\end{thm5}
The curvature of $\Omega'''_{l}$ is given by $R/{\Pi}$, where 
$R =R_{2 {\bar{3}}2{\bar{3}}}-R_{2{\bar{3}}3{\bar{2}}}-
R_{3{\bar{2}}2{\bar{3}}}+R_{3{\bar{2}}3{\bar{2}}}$ and 
$\Pi = 4 <\dot{\mu}_2,\dot{\mu}_2><\dot{\mu}_3,\dot{\mu}_3> - 
2 |<\dot{\mu}_2,\dot{\mu}_3>|^2 - 2 Re(<\dot{\mu}_2,\dot{\mu}_3>)^2$. Here 
we recall that the curvature tensor is 
\begin{center}
$R_{\alpha {\bar{\beta}} \gamma {\bar{\delta}}} =  
(\int_{\Sigma}D(\dot{\mu}_{\alpha} \dot{\bar{\mu}}_{\beta}
\dot{\mu}_{\gamma}\dot{\bar{\mu}}_{\delta}dA) + 
(\int_{\Sigma}D(\dot{\mu}_{\alpha} \dot{\bar{\mu}}_{\delta})
\dot{\mu}_{\gamma}\dot{\bar{\mu}}_{\beta}dA)$
\end{center}
The essential parts on the surface for the curvature estimates are $M_0, M_2$ 
and $M_3$. As in previous sections, we will use rotationally symmetric {\hm}s 
$w_{i}(t)$ to approximate $W_{i}(t)$ in $M_i$, for $i=2, 3$. We still 
denote the {\Hd}s of $w_{i}(t)$ by $\phi_{i}(t)$ and corresponding {\Bd}s 
by $\mu_{i}(t)$ for $i=2, 3$.

Similar to the discussion in model cases one and two, the infinitesimal {\Hd}s 
$\dot{\phi}_{2}(t)$ and $\dot{\phi}_{3}(t)$ are holomorphic. Since 
{\hm}s $w_{2}(t)$ and $w_{3}(t)$ are rotationally symmetric, we can assume the 
infinitesimal {\Hd}s $|\dot{\phi}_{2}|_{M_2} = 1$, 
$|\dot{\phi}_{3}|_{M_3} = 1$, while $|\dot{\phi}_{2}|_{M_0} = \eta_{2}(x,l)$, 
and $|\dot{\phi}_{3}|_{M_0} = \eta_{3} (x,l)$, where we 
consider $M_0$ as a cylinder characterized by $[a,b] \times [0,1]$, and 
$a = a(l) = l^{-1}sin^{-1}(l)$, and 
$b = b(l) = \pi l^{-1}-l^{-1}sin^{-1}(l)$. Here, similar to $\zeta(x,l)$ from 
$\S 4.1$, the functions $\eta_{2}(x,l)$ and $\eta_{3}(x,l)$ satisfy that 
$\eta_{2}(x,l) \leq C_{1}x^{-4}$ for $x \in [a,b]$, and 
$\eta_{3}(x,l) \leq C_{1}(b-x)^{-4}$ for $x \in [a,b]$, and they decay 
exponentially in $[1, +\infty]$. We note that, asymptotically, cylinders 
$M_2$ and $M_3$ lie on different sides of $M_0$. Therefore we can assume that 
$|\dot{\phi}_{2}|_{M_3} = O(l^4)$ and $|\dot{\phi}_{3}|_{M_2} = O(l^4)$.

To prove theorem 1.5, we want to show firstly that $\Pi$ is bounded (away from 
zero). This is easy to see as, asymptotically, both $\dot{\phi}_{2}(t)$ and 
$\dot{\phi}_{3}(t)$ decay exponentially in $M_0$, and hence $\Pi > \Pi|_{M_0}$ 
is bounded from below by some positive constant. Now we just have to show that 
$|R| = O(l)$.

\begin{lem}
$\int_{M_0}D(|\dot{\mu}_2|^2){|\dot{\mu}_3|}^{2}\sigma dxdy = O(l)$
\end{lem}
\begin{proof}

We split the interval $[a,b]$ into three parts: $[a,l^{-1/4}]$, 
$[l^{-1/4}, b-l^{-1/4}]$ and $[b-l^{-1/4},b]$. Recall that 
$D(|\dot{\mu}_2|^2) = \ddot{\mathcal{H}}^M$, where $\ddot{\mathcal{H}}^M$ is the 
holomorphic energy corresponding to the family of {\hm}s $w_{2}(t)$ resulting 
from twisting the curve $\gamma_2$.

Similar to the proof of Lemma 4.3, we compare $\ddot{\mathcal{H}}^M$ with 
$Y(l;x)$, where $Y(l,x) = B_{3}cot(lx) + B_{4}(1-lxcot(lx))$ is the solution 
to $(\Delta-2)Y = 0$ in $[a,b]$ with constants $B_{3}$ and $B_{4}$ 
satisfying that $B_{3} = O(l)$ and $B_{4} = O(l)$. Again, we find that 
$\ddot{\mathcal{H}} + x^{-2}$ is a subsolution to $(\Delta-2)Y = 0$ for 
$x \in [a,b]$. Here the hyperbolic length element in $M_0$ is $lcsc(lx)|dz|$.

In the first interval $[a, l^{-1/4}]$, noticing that 
$|\dot{\phi}_{3}| \le C_{1} (b - l^{-1/4})^{-4} = O(l^4)$, we have 
\begin{eqnarray*}
\int_{a}^{l^{-1/4}}D(|\dot{\mu}_2|^2){|\dot{\mu}_3|}^{2}\sigma dx & \le &
\int_{a}^{l^{-1/4}}(Y(x)-x^{-2})|\dot{\mu}_{3}|^2 \sigma dx \\
& = & \int_{a}^{l^{-1/4}}(Y(x)-x^{-2})|\dot{\phi}_{3}|^{2}l^{-2}sin^{2}(lx)dx \\
& = & o(l)
\end{eqnarray*}
In the second interval $[l^{-1/4}, b-l^{-1/4}]$, we have 
\begin{eqnarray*}
\int_{l^{-1/4}}^{b-l^{-1/4}}D(|\dot{\mu}_2|^2){|\dot{\mu}_3|}^{2}\sigma dx & \le &
\int_{l^{-1/4}}^{b-l^{-1/4}}(Y(x)-x^{-2})|\dot{\mu}_{3}|^2 \sigma dx \\
& \le & \int_{l^{-1/4}}^{b-l^{-1/4}}{\frac{(Y(x)-x^{-2})C_{1}^{2}}{(b-x)^{8}}}
l^{-2}sin^{2}(lx)dx \\
& = & o(l)
\end{eqnarray*}
In the third interval $[b-l^{-1/4}, b]$, noticing that 
$D(|\dot{\mu}_2|^2) \le Y(x) - x^{-2} = O(l)$, we have 
\begin{eqnarray*}
\int_{b-l^{-1/4}}^{b}D(|\dot{\mu}_2|^2){|\dot{\mu}_3|}^{2}\sigma dx & \le &
\int_{b-l^{-1/4}}^{b}(Y(x)-x^{-2})|\dot{\mu}_{3}|^2 \sigma dx \\
& = & \int_{b-l^{-1/4}}^{b}O(l){\frac{C_{1}^{2}}{(b-x)^{8}}}
l^{-2}sin^{2}(lx)dx \\
& = & O(l)
\end{eqnarray*}
Combining these calculation, we complete the proof of Lemma 5.1. 
\end{proof}

To show $|R| = O(l)$, we notice that 
$|R| \le C' \int_{\Sigma}D(|\dot{\mu}_2|^2){|\dot{\mu}_3|}^{2}\sigma dxdy$ for 
some positive constant $C'$, as in the proof of Lemma 3 of {\cite{H}} and 
Lemma 4.2. Hence we need to calculate the integral 
$\int_{\Sigma} D(|\dot{\mu}_2|^2){|\dot{\mu}_3|}^{2} dA$. Since $\gamma_0$ is a 
separating curve, we split this integral into two parts, one for each side of 
$\gamma_0$, and we abuse our notation to denote the side with curve $\gamma_2$ 
by $M_2$, and the other side with curve $\gamma_3$ by $M_3$.

Since $M_2$ is compact, and $|\dot{\phi}_{3}|_{M_2} = O(l^4)$, we have the 
integral 
\begin{eqnarray}
\int_{M_2}D(|\dot{\mu}_2|^2){|\dot{\mu}_3|}^{2}\sigma dxdy = O(l^4) = o(l)
\end{eqnarray}
Also $M_3$ is compact, we have $|\dot{\phi}_{3}|_{M_3} = O(1)$, while 
$D(|\dot{\mu}_2|^2)|_{M_3} = \ddot{\mathcal{H^M}}|_{M_3}$ is comparable to 
$\ddot{\mathcal{H^M}}|_{M_0 \cap M_3}$, which is of the order 
$O(lcsc(lb)) = O(1/b) = O(l)$, recalling that 
$b = b(l) = \pi l^{-1}-l^{-1}sin^{-1}(l)$. Thus 
\begin{eqnarray}
\int_{M_3}D(|\dot{\mu}_2|^2){|\dot{\mu}_3|}^{2}\sigma dxdy = O(l)
\end{eqnarray}
The estimate of Lemma 5.1, (13) and (14) complete the proof of theorem 1.5 
when we consider the family of {\hm}s being rotationally symmetric.

Now we consider the families of {\hm}s which are no longer rotationally 
symmetric. Since $\gamma_0$ separates the surface, the cylinder $M_0$ 
separates the surface into two sides. We denote the side containing 
$\gamma_2$ by $M_2$, and the side containing $\gamma_3$ by $M_3$. As we 
indicated, the infinitesimal {\Hd}s are exponentially small in the thin 
part. Thus we can still assume that $|\dot{\phi}_{2}|_{M_3} = O(l^4)$ 
and $|\dot{\phi}_{3}|_{M_2} = O(l^4)$, hence  the term $\Pi$ in the 
curvature formula is still bounded.

To compute the integral 
$\int_{\Sigma}D(|\dot{\mu}_2|^2) |\dot{\mu}_3|^2 dA$, we still denote 
the holomorphic energy of the family $W_{2}(t)$ by ${\mathcal{H}}(t)$, 
then $D(|\dot{\mu}_2|^2) = \ddot{\mathcal{H}}$ will be dominated, as 
seen in the argument in $\S 3.4$ after (6), by $\ddot{\mathcal{H^M}}$ 
and $Y(l,x) = B_{3}cot(lx) + B_{4}(1-lxcot(lx))$, the solution to 
$(\Delta-2)Y = 0$ in $[a,b]$ with constants $B_{3}$ and $B_{4}$ 
satisfying that $B_{3} = O(l)$ and $B_{4} = O(l)$, when we 
characterize $M_0$ by $[a,b] \times [0,1]$.

Now it is not hard to see that 
\begin{center}
$\int_{\Sigma_2 \cup \Sigma_3}D(|\dot{\mu}_2|^2)
|\dot{\mu}_3|^2 dA = o(l) $ \\
$\int_{[a,b-l^{1/4}] \times [0,1]}D(|\dot{\mu}_2|^2)
|\dot{\mu}_3|^2 dA = o(l) $ 
 \end{center}
While $\int_{[b-l^{1/4},b] \times [0,1]}D(|\dot{\mu}_2|^2)
|\dot{\mu}_3|^2 dA = O(l)$ as 
$Y(l,x) \sim lcot(lx) \sim l$ when $ x \in [b-l^{1/4},b]$. Therefore the 
integral $\int_{\Sigma}D(|\dot{\mu}_2|^2) |\dot{\mu}_3|^2 dA = O(l)$ and 
we complete the proof of theorem 1.5.
\qed
\begin{rem}
If the shrinking curve $\gamma_0 (l)$ is not separating, the twisting 
neighborhoods $M_2$ and $M_3$ will lie on the same component in the limit 
of surfaces with one shrinking curve as $l$ tends to zero. It is not 
hard to see that $K(\Omega'''_l)$ is now bounded away from zero.
\end{rem}
\begin{rem}
Theorem 1.5 can be easily generalized to treat the case when multiple 
closed geodesics are pinching and at least one of them is separating. 
Infinitesimal twists about curves on different sides of a separating and 
shrinking curve will give a family of asymptotically flat {\tp}s. 
\end{rem}
\subsection {Curvature Bounds}
In this subsection, we prove theorem 1.1, i.e., we give curvature bounds at 
any point in {\TS}.
\begin{thm1} 
Let $l$ be the length of the shortest geodesic on closed surface 
$\Sigma$, and $K$ be the {\WP} sectional curvature of 
{\TS} $\mathcal {T}$, assuming $dim_{C}{\mathcal {T}} > 1$, 
there exists a constant $C > 0$ such that 
\begin{center}
$-(Cl)^{-1} \le K \le -Cl$
\end{center}
Moreover, there are {\tp}s with the {\WP} curvatures of 
the orders $O(l)$ and comparable to $l^{-1}$, and hence the {\WP} {\Sc} 
has neither negative upper bound, nor lower bound. 

\end{thm1}

Recall from $\S 3.4$, we showed that the {\WP} holomorphic {\Sc} tends to 
negative infinity at the rate of the order $O(l^{-1})$. We note that, 
asymptotically, the absolute value of the {\Sc} is dominated by diagonal 
terms. One way to see this, assume that $\dot{\nu}_0$ and $\dot{\nu}_1$ 
are two tangent vectors at $\Sigma$ in the tangent space of {\TS}, then, 
as in lemma 4.3 of {\cite {Wp86}}, we have 
$|D(\dot{\nu}_{0}\dot{\bar{\nu}}_1)| \le D(|\dot{\nu}_0|^2)^{1/2}
D(|\dot{\nu}_1|^2)^{1/2}$. 
Applying Schwarz lemma, one finds that the absolute value of the curvature 
of the plane, spanned by $\dot{\nu}_0$ and $\dot{\nu}_1$, is dominated by 
$\int_{\Sigma}D(|\dot{\nu}_0|^2)|\dot{\nu}_1|^2 dA /{\Pi}$, 
where $\Pi = 4 <\dot{\nu}_0,\dot{\nu}_0><\dot{\nu}_1,\dot{\nu}_1> - 
2 |<\dot{\nu}_0,\dot{\nu}_1>|^2 - 2 Re(<\dot{\nu}_0,\dot{\nu}_1>)^2$.

If there is a lower bound for the length of the shortest closed geodesics 
on $\Sigma$, the {\WP} {\Sc} of {\TS} is bounded since all integrals in 
curvature terms are bounded away from negative infinity and zero. Hence, 
to consider large absolute value of the {\Sc}, we can assume one of the 
vectors $\dot{\nu}_0$ and $\dot{\nu}_1$ is corresponding to a deformation 
of the length of a short geodesic on the surface. Recalling that 
the {\Hd} not corresonding to pinching this curve will be exponentially 
small in the thin part of this shrinking curve, so above integral 
$\int_{\Sigma}D(|\dot{\nu}_0|^2)|\dot{\nu}_1|^2 dA$ will be no more than 
the integral $\int_{\Sigma}D(|\dot{\nu}_0|^2)|i\dot{\nu}_0|^2 dA$, which 
we estimated in the proof of theorem 1.2. Therefore the proof of theorem 
1.2 implies the following theorem, which is the lower bound part of 
theorem 1.1, i.e.,
\begin{theorem} 
Let $l$ be the length of the shortest closed geodesic on the surface, 
then there is a positive constant $C$ such that the {\WP} {\Sc} $K$ of 
{\TS} satisfies $K \ge -(Cl)^{-1}$, if the complex dimension of {\TS} is 
great than one.
\end{theorem}
The estimates in proving theorem 1.2 immediately imply the following 
result of Georg Schmuacher:
\begin{cor} (\cite {Sch}) 
The sectional, Ricci and scalar curvature are asymptotically bounded by 
${\sum\limits_{i=1}^{q}}log|t_{i}|$.
\end{cor}
Note that since the {\WP} {\Sc} is negative, therefore this theorem of 
Georg Schmuacher doesnot imply that the curvature is not bounded from below.

To show the upper bound part of theorem 1.1, we notice that the upper 
bounds for asymptotically flat {\tp}s $\Omega_l$ in theorem 1 of 
{\cite{H}}, planes $\Omega'_l$ in theorem 1.3, planes $\Omega''_l$ in 
theorem 1.4, and planes $\Omega'''_l$ in theorem 1.5 are all of the 
order $O(l)$. We want to show $O(l)$ is the right order, in other words, 
we need to show
\begin{theorem}
If $dim_{C}{\mathcal{T}} > 1$, and $l$ is the length of the shortest 
geodesics on the surface, there exists a constant $C > 0$ such that the 
{\WP} {\Sc} $K$ of {\TS} satisfies $K \le -Cl$. 
\end{theorem}
We assume at least one core geodesics is shrinking on the surface. 
Let $(l_1, \theta_1, l_2, \theta_2, ..., \l_{3g-3}, \theta_{3g-3})$ be 
the Fenchel-Nielsen coordinates at a point $\Sigma$ in {\TS}. It 
suffices to prove theorem 5.6 by assuming two tangent vectors 
$\dot{\nu}_0$ and $\dot{\nu}_1$ are infinitesimal {\Bd}s rsulting from 
deformation of lengths of core geodesics or deformation of twisting 
angles or both.

From remark 5.2, if both $\dot{\nu}_0$ and $\dot{\nu}_1$ are resulting 
from deformations of twisting angles about independent core geodesics, 
then this shrinking curve, denoted by $\gamma$ with length $l$, can be
 assumed to be separating (in such case, we require the genus of the 
surface is at least two). Therefore, from the proof of theorem 1.5, 
remark 5.2 and remark 5.3, we have 
$K(Span(\dot{\nu}_0,\dot{\nu}_1))= O(l)$.

Now we assume at least one of $\dot{\nu}_0$ and $\dot{\nu}_1$ is 
resulting from a deformation of the length of a core geodesic on the 
surface, with or without twisting about this geodesic. Since operator 
$D = -2(\Delta - 2)^{-1}$ is self-adjoint, we assume $\dot{\nu}_0$ is 
resulting from a deformation of the length of a core geodesic 
$\gamma_0$. As in the argument before theorem 5.4, we need to 
estimate the integral 
$\int_{\Sigma}D(|\dot{\nu}_0|^2)|\dot{\nu}_1|^2 dA$ as it will 
dominate the absolute value of $R =  R_{0 {\bar{1}} 0 {\bar{1}}} - 
R_{0 {\bar{1}} 1 {\bar{0}}} - R_{1 {\bar{0}} 0 {\bar{1}}} + 
R_{1 {\bar{0}} 1 {\bar{0}}}$.

Let $M_0$ be the pinching neighborhood of the shrinking geodesic 
$\gamma_0$ with $l(\gamma_0)=l$, as before, we can characterize 
$M_0$ as $[a,b] \times [0,1]$, where 
$a = a(l) = l^{-1}sin^{-1}(l)$, and 
$b = b(l) = \pi l^{-1}-l^{-1}sin^{-1}(l)$. Therefore we find that 
$D(|\dot{\nu}_0|^2)$ is of the order of $O(lcot(lb)) = O(l)$ in 
$\Sigma \backslash M_0$. Hence,
\begin{eqnarray*}
\int_{\Sigma}D(|\dot{\nu}_0|^2)|\dot{\nu}_1|^2 dA & = & 
\int_{M_0}D(|\dot{\nu}_0|^2)|\dot{\nu}_1|^2 dA + 
\int_{\Sigma \backslash M_0}D(|\dot{\nu}_0|^2)|\dot{\nu}_1|^2 dA \\
& = & \int_{M_0}D(|\dot{\nu}_0|^2)|\dot{\nu}_1|^2 dA + 
O(l)\int_{\Sigma \backslash M_0}|\dot{\nu}_1|^2 dA 
\end{eqnarray*}
Note that $\int_{\Sigma \backslash M_0}|\dot{\nu}_0|^2 dA = O(1)$, 
we have the absolute value of the curvature is given by
\begin{eqnarray*}
|R|/{\Pi} & \le & 
{\frac{C'\int_{\Sigma}D(|\dot{\nu}_0|^2)|\dot{\nu}_1|^2 dA}{\Pi}} \\ 
& \le & C'{\frac{\int_{M_0}D(|\dot{\nu}_0|^2)|\dot{\nu}_1|^2 dA}
{\Pi}} + O(l){\frac{\int_{\Sigma \backslash M_0}|\dot{\nu}_1|^2 
dA}{\Pi}} \\
& = & C'({\frac{\int_{M_0}D(|\dot{\nu}_0|^2)|\dot{\nu}_1|^2 dA}{\Pi}}) 
+ O(l)
\end{eqnarray*}
We notice that the integral 
$\int_{M_0}D(|\dot{\nu}_0|^2)|\dot{\nu}_1|^2 dA$ is positive, and will 
be smaller when the infinitesimal {\Bd}s $\dot{\nu}_1$ is resulting 
from a deformation of the length of a core geodesic, independent of 
$\gamma_0$, on the surface, than otherwise. To see this, 
we recall that $|\dot{\phi}_1|$, where $\dot{\phi}_1$ is the 
infinitesimal {\Hd} corresponding to $\dot{\nu}_1$, decays 
exponentially away from the curve where the deformation occurs. 
Therefore $|\dot{\phi}_1|$ is smaller in $M_0$ when it is 
actually resulting from a deformation of the length of a core geodesic 
$\gamma_1$ than the case of otherwise. Now the proofs of theorem 1 of 
{\cite{H}}, theorem 1.3, and theorem 1.4 imply that 
\begin{eqnarray*}
|R|/{\Pi} & \le & C'(
{\frac{\int_{M_0}D(|\dot{\nu}_0|^2)|\dot{\nu}_1|^2 dA}{\Pi}}) + O(l) \\ 
& = & O(l) + O(l) = O(l)
\end{eqnarray*}
This completes the proof of theorem 5.6. Theorem 1.1 immediately 
follows from theorem 5.4 and theorem 5.6.

\noindent
Zheng Huang\\
Department of Mathematics \\
University of Michigan\\
Ann Arbor, MI 48109\\
email address: zhengh@umich.edu
\end{document}